\documentclass[10pt,a4paper]{article}
\usepackage{amsfonts,amsmath,mathrsfs,amssymb}
\usepackage{tikz}
\usepackage{verbatim}
\usepackage{array}
\newtheorem{theorem}{Theorem}[section]

\newtheorem{prop}[theorem]{Proposition}
\newtheorem{defin}{Definition}[section]
\newtheorem{examp}{Example}[section]
\newtheorem{remark}{Remark}[section]
\newtheorem{cor}[theorem]{Corollary}

\usepackage{pxfonts}
\usepackage{verbatim}
\usepackage{array}
\usepackage{tabularx}
\usepackage{subfig}

\def\R{ \mathbb{R}}
\let\tilde\widetilde
\let\hat\widehat
\let\bar\overline

\def\op{\overline{p}}
\def\ps{{\op}^*}
\def\RC{ {R_*}}
\def\Rb{\bar{R}}
\def\Ob{\bar{O}}
\def\OC{ {O_*}}

 \def\SS{\mathbb S}
 \def\M{\mathbb M}
 \def\H{\mathbb H}
  \def\E{\mathbb E}
 \def\DS{{d}\mathbb{S}}

\definecolor{pink}{rgb}{1, .509, .8}
\definecolor{sky}{rgb}{.4, .914, 1}

\begin{document}

\title {Coning, Symmetry and Spherical Frameworks}

\author{Bernd Schulze\footnote{Supported by the DFG Research Unit 565 `Polyhedral Surfaces'.}\\ Fields Institute for Research in Mathematical Sciences\\ 222 College Street\\ Toronto, ON M5T-3J1, Canada\\
and\\ Walter Whiteley\footnote{Supported by a grant from NSERC (Canada).}\\
 Department of Mathematics and Statistics\\ York University\\ 4700 Keele Street\\ Toronto, ON M3J-1P3, Canada
}

\maketitle

\begin{abstract}
In this paper, we combine separate works on (a) the transfer of infinitesimal rigidity results from an Euclidean space to the next higher dimension by coning \cite{WWcones}, (b) the further transfer of these results to spherical space via associated rigidity matrices \cite{SalWW}, and (c) the prediction of finite motions from symmetric infinitesimal motions at regular points of the symmetry-derived orbit rigidity matrix \cite{BSWWorbit}.   Each of these techniques is reworked and simplified to apply across several metrics, including the Minkowskian metric $\M^{d}$ and the hyperbolic metric $\H^{d}$.

This leads to a set of new results transferring infinitesimal and finite motions associated with corresponding symmetric frameworks among $\E^{d}$, cones in $E^{d+1}$, $\SS^{d}$, $\M^{d}$, and $\H^{d}$.  We also consider the further extensions associated with the  other Cayley-Klein geometries overlaid on the shared underlying projective geometry.

\end{abstract}

\section{Introduction}
It is a difficult problem to predict when a framework $(G,p)$ in Euclidean space $\E^{d}$ is flexible.  It is even less common to predict that a framework $(G,p)$ is flexible in another metric, such as the spherical space $\SS^{d}$, hyperbolic space $\H^{d}$, or even Minkowskian space $\M^{d}$.  There are a few well known examples, such as the Bricard octahedra, which have been shown to be flexible in each of these spaces, but with a separate proof for each space \cite{bricard,alex,Stachel,Stachelhyp}.  Such examples raise the possibility that there are transfer principles which would bring some flexibility results from one space to another.   In this paper, we present one such transfer principle.

Coning is an established tool for transferring results about infinitesimal rigidity of a framework from $\E^{d}$ to $\E^{d+1}$ \cite{WWcones}.  Working within the cone, with the original framework placed on the hyperplane $x_{d+1}=1$, then pulling points back and forth towards the cone vertex, we can also transfer infinitesimal rigidity results from $\E^{d}$ through the cone to the spherical  space $\SS^{d}$ within $\E^{d+1}$, with the cone point at the center of the sphere  \cite{ SalWW}.   As a more subtle consequence of the equivalence of infinitesimal flexibility and finite flexibility at regular configurations (configurations where the rigidity matrix indexed by the edges and vertices of the framework has maximal rank), this process also transfers finite flexes between  $\E^{d}$,  the cone in $\E^{d+1}$,  and the sphere $\SS^{d}$.

Several recent papers have examined when symmetry of a framework induces unexpected finite flexes in $\E^{d}$, most recently using a symmetric analog of the rigidity matrix: the orbit matrix \cite{BSWWorbit}.   This matrix has rows indexed by orbits of edges under the symmetry group, and columns indexed by orbits of vertices under the group.  Again, in this setting, fully symmetric infinitesimal flexes and finite flexes are equivalent at symmetry regular configurations (configurations where the orbit matrix has maximal rank).

It is a natural question to ask when this symmetry analysis is preserved by coning, so that a symmetry induced finite flex in $\E^{d}$
transfers to the cone in $\E^{d+1}$  and to the sphere $\SS^{d}$. Our starting point in \S3 is to confirm that this transfer works when the cone point is placed  `above' the central point of the symmetry group (point group).   In the process, we give an alternate proof of the transfer for the identity group, using direct matrix manipulations on the rigidity matrix, and then extend the analysis to the orbit matrix.

Previous work on rigidity in non-Euclidean geometry has included results on the transfer to hyperbolic space \cite{SalWW}.   We fit this into our scheme by (a) embedding the framework $(G,p)$ in $\E^{d}$ on the hyperplane $x_{d+1}=1$ in Minkowskian space $\M^{d+1}$ and (b) coning the framework to the origin, in a way that transfers all the infinitesimal rigidity properties.  This works for an arbitrary symmetry group in $\E^{d}$, carrying over to symmetries in $\M^{d+1}$, and preserving the rank of the corresponding orbit matrix.  Now, when we push the vertices onto the `unit sphere', we have a framework in hyperbolic space $\H^{d}$.  This preserves all the key properties including taking symmetry regular configurations in $\E^{d}$ onto symmetry regular configurations of the cone in  $\M^{d+1}$, and finite flexes at regular points to finite flexes at regular points.

These results suggest two directions for extensions.  One is: What other metrics can we transfer to?  We recall that infinitesimal rigidity (but not symmetry) is a projective invariant \cite{CW}, and can be expressed in projective terms.  It  is then natural to consider other metrics layered onto projective geometry - the Cayley-Klein geometries \cite{SalWW,CayleyKlein,richter}.  In \S\ref{sec:projective} we briefly outline how the transfer process generalizes to families of these geometries, provided that the metric is associated with a quadratic form with a signature of all positive and negative entries, or with a sphere in such a space.   Exactly how far this transfer stretches is a topic for further explorations.

A second direction for extension is the big question: when do finite motions transfer from frameworks in one space to frameworks in another space?  It is a hard problem to decide whether a specific framework $(G,p)$ which is first-order flexible in a special position is rigid or flexible.  It is then at least as hard to decide whether flexibility would transfer to another metric with, for example, the same projective coordinates.

The results in this paper indicate that this transfer can be accomplished when we have an algebraic variety of configurations $\cal{A}$, within which the framework moves for regular points of the variety.  The key case here is the variety of $S$-symmetric configurations for a group $S$, where symmetric infinitesimal motions and the symmetry regular configurations are detected by the orbit matrix.  Further, the algebraic variety is preserved in the transfer, as is the rank of an associated matrix which tests for the infinitesimal motions tangent to the variety.   There are some other varieties which are candidates for such transfer.  One sample  is the varieties associated with `flatness' - sets of points in the configuration lying in linear subspaces, which are known to generate finite motions for some classes of complete bipartite graphs.   Again, the range of extensions of these methods is an inviting question for further exploration.


\subsection{Overview Tables}

We summarize the key results of this paper both for frameworks without symmetry, and frameworks with symmetry in several tables.

\begin{table}[h!]
\begin{tabularx}{.4 \textwidth}{X}
\[ \begin{array}{|X||c||X|}
\hline
{\textrm{Framework $(G,p)$ in $\E^{d}$}}  &   {\textrm{$p=\pi(q)$}}  & {\textrm{Coned framework $(G*o,q)$ in $\E^{d+1}$}} \\
 \hline
 \hline
\textrm{dim (trivial infinitesimal motions of $(G,p)$)$=m$ } &  \Leftrightarrow & \textrm{dim (trivial infinitesimal motions of $(G*o,q)$) $=m+(d+1)$} \\
\hline
\textrm{$(G,p)$ has a non-trivial infinitesimal motion } &  \Leftrightarrow & \textrm{$(G*o,q)$ has a non-trivial infinitesimal motion}\\
\hline
\textrm{$(G,p)$ has a non-trivial self-stress} &   \Leftrightarrow  & \textrm{$(G*o,q)$ has a non-trivial self-stress} \\
\hline
\textrm{$p$ is a regular point of $G$ } &   \Leftrightarrow  & \textrm{$q$ is a regular point of $G*o$} \\
\hline
\textrm{$p$ is a regular point with a finite motion } &   \Leftrightarrow  & \textrm{$q$ is regular point a with a finite motion} \\
\hline
\end{array}  \]
\end{tabularx}
\label{table:outline}
  \caption{Summary of results about the transfer of flexibility from the framework $(G,p)$ in the Euclidean space $\E^{d}$ to the coned framework $(G*o,q)$ in $\E^{d+1}$ (see Theorems \ref{thm:ConingGeneral}, \ref{thmtranconwithsym}, \ref{thm:regptstransfer}, and \ref{thm:coning flexes}).}
\end{table}

\begin{table}[h!]
\begin{tabularx}{.4  \textwidth}{X}
\[ \begin{array}{|X||c||X|}
\hline
{\textrm{Framework $(G,p)$ with symmetry group $S$ in $\E^{d}$}}  &   {\textrm{$p=\pi(q)$}}  & {\textrm{Coned framework $(G*o,q)$ with symmetry group $S^{*}$ in $\E^{d+1}$}} \\
 \hline
 \hline
\textrm{$(G,p)$ has a non-trivial $S$-symmetric infinitesimal motion} &  \Leftrightarrow & \textrm{$(G*o,q)$  has a non-trivial $S^{*}$-symmetric infinitesimal motion} \\
\hline
\textrm{$(G,p)$ has a non-trivial $S$-symmetric self stress} &   \Leftrightarrow  & \textrm{$(G*o,q)$ has a non-trivial $S^{*}$-symmetric self stress} \\
\hline
\textrm{$p$ is an $S$-regular point of $G$ } &   \Leftrightarrow  & \textrm{$q$ is an $S^{*}$-regular point of $G*o$} \\
\hline
\textrm{$p$ is an $S$-regular point with symmetry-preserving finite motion } &   \Leftrightarrow  & \textrm{$q$ is an $S^{*}$-regular point with symmetry-preserving finite motion} \\
\hline
\end{array}  \]
\end{tabularx}
\label{table:outline2}
  \caption{Summary of results about the transfer of $S$-symmetric flexibility from the framework $(G,p)$ with symmetry group $S$ in the Euclidean space $\E^{d}$ to the coned framework $(G*o,q)$ with symmetry group $S^*$ in $\E^{d+1}$ (see Theorems \ref{thmtranconwithsym}, \ref{thm:regptstransfer}, and \ref{thm:coning flexes}). }
\end{table}

\begin{table}[h!]
\begin{tabularx}{.4  \textwidth}{X}
\[ \begin{array}{|X||c||X|}
\hline
{\textrm{Framework $(G,p)$ with symmetry group $S$ in $\E^{d}$}}  &   {\textrm{$p=\pi(q)$}}  & {\textrm{Framework $(G,q)$ with symmetry group $S^{*}$ in $\SS^{d}$}} \\
 \hline
 \hline
\textrm{$(G,p)$ has a non-trivial $S$-symmetric infinitesimal motion} &  \Leftrightarrow & \textrm{$(G,q)$  has a non-trivial $S^{*}$-symmetric infinitesimal motion} \\
\hline
\textrm{$(G,p)$ has a non-trivial $S$-symmetric self stress} &   \Leftrightarrow  & \textrm{$(G,q)$ has a non-trivial $S^{*}$-symmetric self stress} \\
\hline
\textrm{$p$ is an $S$-regular point of $G$ } &   \Leftrightarrow  & \textrm{$q$ is an $S^{*}$-regular point of $G$} \\
\hline
\textrm{$p$ is an $S$-regular point with symmetry-preserving finite motion } &   \Leftrightarrow  & \textrm{$q$ is an $S^{*}$-regular point with symmetry-preserving finite motion} \\
\hline
\end{array}  \]
\end{tabularx}
\label{table:outline3}
  \caption{Summary of results about the transfer of $S$-symmetric flexibility from the framework $(G,p)$ with symmetry group $S$ in the Euclidean space $\E^{d}$ to the framework $(G,q)$ with symmetry group $S^*$ in the spherical space $\SS^{d}$. (see Theorems \ref{thm:hemisphere}, \ref{thm:regptstransfer}, \ref{thm:coning flexes}, and Corollary \ref{thm:wholesphere}).}
\end{table}

\section{Background }

\subsection{Rigidity of frameworks in $\E^{d}$}

The following is standard material, following sources such as \cite{W1,BS2}. We define a \emph{framework}  in $\mathbb{E}^{d}$ to be a pair $(G,p)$, where $G$ is a finite simple graph with vertex set $V(G)$ and edge set $E(G)$, and  $p:V(G)\to \mathbb{E}^d$ is an embedding of the vertices of $G$ in Euclidean $d$-space.  We also say that $(G,p)$ is a $d$-dimensional \emph{realization} of the \emph{underlying graph} $G$. We often identify the function $p$ with a row vector in $\mathbb{E}^{d|V(G)|}$ (by using some fixed order
on the vertices in $V(G)$), in which case we refer to $p$ as a \emph{configuration} of $|V(G)|$ points in $\mathbb{E}^{d}$.
For $i\in V(G)$, we say that $p(i)$ is the \emph{joint} of $(G,p)$ corresponding to $i$, and for $e=\{i,j\}\in E(G)$, we say that $\{p(i),p(j)\}$ is the \emph{bar} of $(G,p)$ corresponding to $e$. Throughout the paper, we denote $n:=|V(G)|$.

For a fixed ordering of the edges of a graph $G$, we define the \emph{edge function} $f_{G}:\mathbb{E}^{dn}\to \mathbb{R}^{|E(G)|}$  by
\begin{displaymath}
f_{G}\big(p_{1},\ldots,p_{n}\big)=\big(\ldots, \|p_{i}-p_{j}\|^2,\ldots\big)\textrm{, }
\end{displaymath}
where $\{i,j\}\in E(G)$, $p_{i}:=p(i)\in \mathbb{E}^{d}$ for all $i\in V(G)$, and $\| \cdot \|$ denotes the Euclidean norm in
$\mathbb{E}^{d}$  \cite{asiroth}.

If $(G,p)$ is a $d$-dimensional framework, then $f_{G}^{-1}\big(f_{G}(p)\big)$ is the set of all configurations $q$ of $n$ points in $\mathbb{E}^{d}$ with the property that corresponding bars of the frameworks $(G,p)$ and $(G,q)$ have the same length. In particular, we clearly have
$f^{-1}_{K_{n}}\big(f_{K_{n}}(p)\big)\subseteq f_{G}^{-1}\big(f_{G}(p)\big)$, where $K_{n}$ is the complete graph on $V(G)$.

An analytic path $x:[0,1]\to \mathbb{E}^{dn}$ is called a \emph{finite motion} of $(G,p)$ if $x(0)=p$ and $x(t)\in f_{G}^{-1}\big(f_{G}(p)\big)$ for all
$t\in [0,1]$. Further, $x$ is called a \emph{finite rigid motion} (or \emph{trivial finite motion}) if $x(t)\in f_{K_{n}}^{-1}\big(f_{K_{n}}(p)\big)$ for all $t\in [0,1]$, and $x$ is called a \emph{finite flex} (or \emph{non-trivial finite motion}) of $(G,p)$ if $x(t)\notin f_{K_{n}}^{-1}\big(f_{K_{n}}(p)\big)$ for all $t\in (0,1]$.

We say that $(G,p)$ is \emph{rigid} if every finite motion of $(G,p)$ is trivial; otherwise $(G,p)$ is called \emph{flexible}. It is a well established fact that the existence of an analytic finite flex is equivalent to the existence of a continuous finite flex, and in turn to a converging sequence of non-congruent configurations \cite{asiroth}.

Given a framework $(G,p)$ in $\E^{d}$, the
\emph{rigidity matrix} of $(G,p)$ is the $|E(G)| \times dn$ matrix $\mathbf{R}(G,p)=\frac{1}{2}df_G(p)$, where $df_G(p)$ denotes the Jacobian matrix of the edge function $f_G$, evaluated at the point $p$. We have
\begin{displaymath} \mathbf{R}(G,p)=\bordermatrix{& & & & i & & & & j & & & \cr & & & &  & & \vdots & &  & & &
\cr \{i,j\} & 0 & \ldots &  0 & (p_{i}-p_{j}) & 0 & \ldots & 0 & (p_{j}-p_{i}) &  0 &  \ldots&  0 \cr & & & &  & & \vdots & &  & & &
}
\textrm{.}\end{displaymath}

An \emph{infinitesimal motion} of a framework $(G,p)$ in $\mathbb{E}^d$ is a function $u: V(G)\to \mathbb{E}^{d}$ such that
\begin{equation}
\label{infinmotioneq}
(p_i-p_j)\cdot (u_i-u_j)=0 \quad\textrm{ for all } \{i,j\} \in E(G)\textrm{,}\end{equation}
where $u_i$ denotes the vector $u(i)$ for each $i$. Note that the kernel of the rigidity matrix $\mathbf{R}(G,p)$ is the space of all infinitesimal motions of $(G,p)$.

An infinitesimal motion $u$ of $(G,p)$ is an \emph{infinitesimal rigid motion} (or \emph{trivial infinitesimal motion}) if there exists a skew-symmetric matrix $S$ (a rotation) and a vector $t$ (a translation) such that $u_i=Sp_i+t$ for all $i\in V(G)$. Otherwise $u$ is an \emph{infinitesimal flex} (or \emph{non-trivial infinitesimal motion}) of $(G,p)$.

Note that if the joints of $(G,p)$ span all of $\mathbb{E}^d$ (in an affine sense), then the kernel of the rigidity matrix $\mathbf{R}(K_n,p)$, where $K_n$ is the complete graph on the vertices of $G$, is the space of infinitesimal rigid motions of $(G,p)$. It is well known that this space is of dimension $\binom{d+1}{2}$.

We say that $(G,p)$ is \emph{infinitesimally rigid} if every infinitesimal motion of $(G,p)$ is an infinitesimal rigid motion. Otherwise $(G,p)$ is said to be \emph{infinitesimally flexible} \cite{W1}.

Clearly, if the joints of $(G,p)$ affinely span all of $\mathbb{E}^d$, then $\textrm{nullity }\big(\mathbf{R}(G,p)\big)\geq \binom{d+1}{2}$, and $(G,p)$ is infinitesimally rigid if and only if $\textrm{nullity } \big(\mathbf{R}(G,p)\big)=\binom{d+1}{2}$ or equivalently, $\textrm{rank }\big(\mathbf{R}(G,p)\big)=d n - \binom{d+1}{2}$.

An infinitesimally rigid framework is always rigid. The converse, however, does not hold in general.

A \emph{self-stress} of a framework $(G,p)$ is a function  $\omega:E(G)\to \mathbb{E}$ such that at each joint $p_i$ of $(G,p)$ we have
\begin{displaymath}
\sum_{j :\{i,j\}\in E(G)}\omega_{ij}(p_{i}-p_{j})=0 \textrm{,}
\end{displaymath}
where $\omega_{ij}$ denotes $\omega(\{i,j\})$ for all $\{i,j\}\in E(G)$. Note that if we identify a self-stress $\omega$ with a row vector in $\mathbb{R}^{|E(G)|}$ (by using the order on $E(G)$), then we have $\omega\mathbf{R}(G,p)=0$.
In structural engineering, the self-stresses are also called \emph{equilibrium stresses} as they record  tensions and compressions in the bars balancing at each vertex.

If $(G,p)$ has a non-zero self-stress, then $(G,p)$ is said to be \emph{dependent} (since in this case there exists a linear dependency among the row vectors of $\mathbf{R}(G,p)$). Otherwise, $(G,p)$ is said to be \emph{independent}.  A framework which is both independent and infinitesimally rigid is called \emph{isostatic} \cite{CW, W1}.

A configuration $p$ of $n$ points in $\mathbb{E}^d$ is called a \emph{regular point of the graph $G$}, if  $\textrm{\rm rank\,}\big(\mathbf{R}(G,p)\big)\geq\textrm{\rm rank\,} \big(\mathbf{R}(G,q)\big)$ for all $q\in \mathbb{E}^{dn}$.
A framework $(G,p)$ is said to be \emph{regular} if $p$ is a regular point of $G$.

It follows from this definition that the set of all regular realizations of a graph $G$ in $\mathbb{E}^d$ forms a dense open subset of all possible realizations of $G$ in $\mathbb{E}^d$ \cite{gss}. Moreover, note that the infinitesimal rigidity of a regular realization of $G$ depends only on the underlying graph $G$ and not on the particular realization \cite{gss}.

Asimov and Roth showed in \cite{asiroth} that for regular frameworks, infinitesimal rigidity and rigidity are  equivalent. This result of Asimov and Roth provides a key tool for detecting finite motions in frameworks.
An extension of the theorem of Asimov and Roth to symmetric frameworks has recently been established in \cite{BS6}. We will formulate this result in the next section.

     \subsection{Symmetric frameworks in $\E^{d}$}
\label{sec:sym}

Let $G$ be a graph and let $\textrm{Aut}(G)$ denote the automorphism group of $G$. A \emph{symmetry operation} of a framework $(G,p)$ in $\mathbb{E}^{d}$ is an isometry $x$ of $\mathbb{E}^{d}$ such that for some $\alpha\in \textrm{Aut}(G)$, we have
$x(p_i)=p_{\alpha(i)}$ for all $i\in V(G)$ \cite{Hall, BS2, BS1}.

The set of all symmetry operations of a framework $(G,p)$ forms a group under composition, called the \emph{point group} of $(G,p)$ \cite{bishop, Hall, BS1}. Since translating a framework does not change its rigidity properties, we may assume wlog that the point group of any framework in this paper is a \emph{symmetry group} (with the origin fixed), i.e., a subgroup of the orthogonal group $O(\mathbb{E}^{d})$ \cite{BS6, BS2, BS1}.

We use the Schoenflies notation for the symmetry operations and symmetry groups considered in this paper, as this is one of the standard notations in the literature about symmetric structures (see \cite{bishop, FG4, Hall, KG1, KG2, BS6, BS2}, for example). In this notation, the identity transformation is denoted by $Id$, a rotation  about a $(d-2)$-dimensional subspace of $\mathbb{E}^d$ by an angle of $\frac{2\pi}{m}$ is denoted by $C_m$, and a reflection in a $(d-1)$-dimensional subspace of $\mathbb{E}^d$ is denoted by $s$.

While the general methods and results of this paper apply to all symmetry groups, we will only analyze three different types of groups in our examples.  In the Schoenflies notation, they are denoted by $\mathcal{C}_{m}$, $\mathcal{C}_{s}$, and $\mathcal{C}_{mv}$.  For any dimension $d$, $\mathcal{C}_{m}$ is a symmetry group generated by an $m$-fold rotation $C_m$, and $\mathcal{C}_{s}$ is a symmetry group consisting of the identity $Id$ and a reflection $s$. The only other possible type of symmetry group in dimension $2$ is the group $\mathcal{C}_{mv}$ which is a dihedral group generated by a pair $\{C_m, s\}$. In dimension $3$, $\mathcal{C}_{mv}$ denotes any symmetry group that is generated by a rotation $C_m$ and a reflection $s$ whose corresponding mirror contains the rotational axis of $C_m$. For further information about the Schoenflies notation we refer the reader to \cite{bishop, Hall, BS1}.

\begin{figure}[htp]
\begin{center}
\begin{tikzpicture}[very thick,scale=1]
\tikzstyle{every node}=[circle, draw=black, fill=white, inner sep=0pt, minimum width=5pt];
        \path (1.5,-1) node (p2) [label = below left: $p_{2}$] {} ;
       \path (0.6,0.5) node (p1) [label = left: $p_{1}$] {} ;
   \path (2.4,0.5) node (p3) [label = right: $p_{3}$] {} ;
   \path (1.5,1.1) node (p4) [label = above left: $p_{4}$] {} ;
   \draw (p1) -- (p4);
      \draw (p3) -- (p4);
     \draw (p2) -- (p3);
      \draw (p2) -- (p1);
        \draw [dashed, thin] (1.5,-1.6) -- (1.5,1.6);
      \node [draw=white, fill=white] (a) at (1.5,-2.1) {(a)};
    \end{tikzpicture}
    \hspace{2cm}
        \begin{tikzpicture}[very thick,scale=1]
\tikzstyle{every node}=[circle, draw=black, fill=white, inner sep=0pt, minimum width=5pt];
    \path (-0.7,0.8) node (p1) [label = above left: $p_1$] {} ;
    \path (0.7,0.8) node (p4) [label = above right: $p_4$] {} ;
    \path (-1.6,-0.8) node (p2) [label = below left: $p_2$] {} ;
     \path (1.6,-0.8) node (p3) [label = below right: $p_3$] {} ;
      \draw (p1) -- (p4);
    \draw (p1) -- (p2);
    \draw (p3) -- (p4);
    \draw (p2) -- (p3);
     \draw [dashed, thin] (0,-1.6) -- (0,1.6);
      \node [draw=white, fill=white] (b) at (0,-2.1) {(b)};
        \end{tikzpicture}
\end{center}
\vspace{-0.3cm}
\caption{\emph{$2$-dimensional realizations of $K_{2,2}$ in $\mathscr{R}_{(K_{2,2},\mathcal{C}_s)}$ of different types: the framework in \emph{(a)} is of type
$\Phi_{a}$, where $\Phi_{a}: \mathcal{C}_{s} \to \textrm{Aut}(K_{2,2})$ is the homomorphism defined by $\Phi_{a}(s)=
(1 \, 3)(2)(4)$ and the framework in \emph{(b)} is of type
$\Phi_{b}$, where $\Phi_{b}: \mathcal{C}_{s} \to \textrm{Aut}(K_{2,2})$ is the homomorphism defined by $\Phi_{b}(s)=
(1 \, 4)(2\, 3)$.}}
\label{K33types}
\end{figure}
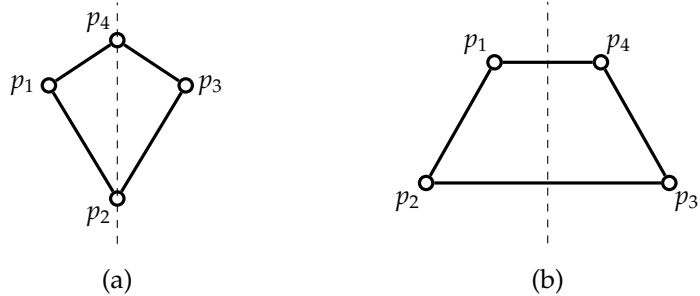

Given a symmetry group $S$ in dimension $d$ and a graph $G$, we let $\mathscr{R}_{(G,S)}$ denote the set of all $d$-dimensional realizations of $G$ whose point group is either equal to $S$ or contains $S$ as a subgroup \cite{BS2, BS1}. In other words, the set $\mathscr{R}_{(G,S)}$ consists of all realizations $(G,p)$ of $G$ for which there exists a map $\Phi:S\to \textrm{Aut}(G)$ so that
\begin{equation}\label{class} x(p_i)=p_{\Phi(x)(i)}\textrm{ for all } i\in V(G)\textrm{ and all } x\in S\textrm{.}\end{equation}
A framework $(G,p)\in \mathscr{R}_{(G,S)}$ satisfying the equations in (\ref{class}) for the map $\Phi:S\to \textrm{Aut}(G)$ is said to be \emph{of type $\Phi$}, and the set of all realizations in $\mathscr{R}_{(G,S)}$ which are of type $\Phi$ is denoted by $\mathscr{R}_{(G,S,\Phi)}$ (see again \cite{BS2, BS1} and Figure \ref{K33types}).

Since we assume that the map $p$ of any framework $(G,p)\in \mathscr{R}_{(G,S,\Phi)}$ is always injective (i.e., $p_i\neq p_j$ if $i\neq j$), it follows that $(G,p)$ is of a unique type $\Phi$, and $\Phi$ is necessarily also a homomorphism (see \cite{BS1} for details). This allows us (with a slight abuse of notation) to use the terms $p_{x(i)}$ and  $p_{\Phi(x)(i)}$ interchangeably, where $i\in V(G)$ and $x\in S$. In general, if the type $\Phi$ is clear from the context, we often simply write $x(i)$ instead of $\Phi(x)(i)$.

Let $(G,p)\in \mathscr{R}_{(G,S,\Phi)}$ and let  $x$ be a symmetry operation in $S$. Then the joint $p_i$ of $(G,p)$ is said to be \emph{fixed} by $x$ if $\Phi(x)(i)=i$, or equivalently (since $p$ is assumed to be injective), $x(p_i)=p_i$.

Let the \emph{symmetry element} corresponding to  $x$ be the linear subspace $F_{x}$ of $\mathbb{E}^{d}$ which consists of all points $a\in\mathbb{E}^{d}$ with $x(a)=a$. Then the joint $p_i$ of any framework $(G,p)$ in $\mathscr{R}_{(G,S,\Phi)}$ must lie in the linear subspace $$U(p_i)=\bigcap_{x\in S: x(p_i)=p_i}F_x\textrm{.}$$

The joint $p_1$ of the framework $(K_{2,2},p)\in \mathscr{R}_{(K_{2,2},\mathcal{C}_s,\Phi_a)}$ depicted in Figure \ref{K33types} (a), for example, is fixed by the identity $Id\in \mathcal{C}_s$, but not by the reflection $s\in \mathcal{C}_s$, so that $U(p_1)=F_{Id}=\mathbb{E}^2$. The joint $p_2$ of $(K_{2,2},p)$, however, is fixed by both the identity $Id$ and the reflection $s$ in  $\mathcal{C}_s$, so that $U(p_2)=F_{Id}\cap F_s=F_s$. In other words, $U(p_2)$ is the mirror line corresponding to $s$.

Note that if we choose a set of representatives $\mathscr{O}_v=\{1,\ldots, k\}$ for the orbits  $S(i)=\{\Phi(x)(i)|\,x\in S\}$ of vertices of $G$, then the positions of \emph{all} joints of  $(G,p)\in \mathscr{R}_{(G,S,\Phi)}$ are uniquely determined by the positions of the joints $p_1,\ldots,p_k$ and the symmetry constraints imposed by $S$ and $\Phi$. Thus, any framework in $\mathscr{R}_{(G,S,\Phi)}$ may be constructed by first choosing positions $p_i\in U(p_{i})$ for each $i=1,\ldots, k$, and then letting $S$ and $\Phi$ determine the positions of the remaining joints.

Let $\mathscr{P}_{(G,S,\Phi)}$ be the subspace of configurations in $\mathbb{E}^{dn}$ which satisfy the equations in (\ref{class}). Further, let $\tilde{f}_{G}: \mathscr{P}_{(G,S,\Phi)}\to\mathbb{E}^{|E(G)|}$ denote the restriction of the edge function $f_{G}$ to $\mathscr{P}_{(G,S,\Phi)}$. A configuration $p\in\mathscr{P}_{(G,S,\Phi)}$ is called  an \emph{$S$-regular point of $G$} if
$\textrm{\rm rank\,}\big(d\tilde{f}_{G}(p)\big)={\max }\{\textrm{\rm rank\,}\big(d\tilde{f}_{G}(q)\big)|\, q\in \mathscr{P}_{(G,S,\Phi)}\}$. A framework $(G,p)\in \mathscr{R}_{(G,S,\Phi)}$ is called $S$-regular if the configuration $p$ is an $S$-regular point of $G$.

To formulate the symmetric version of the theorem of Asimov and Roth, we need the notion of an $S$-symmetric infinitesimal motion which we define next. In general, as shown in \cite{BS6,BSWWorbit}, an analysis of the infinitesimal motions and stresses of a symmetric framework which exhibit the full symmetry of the framework can give important insight into the rigidity properties of the framework.

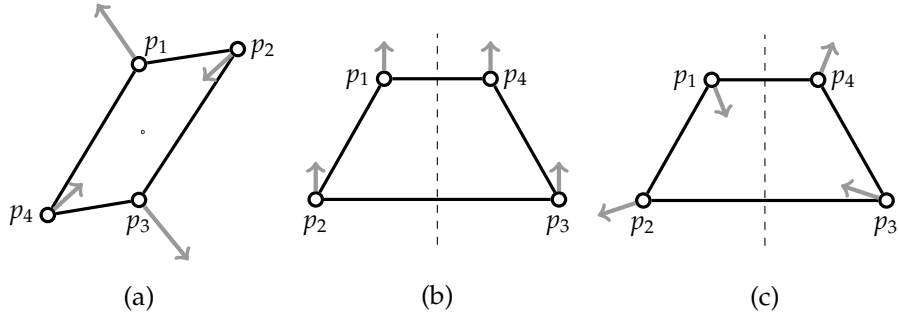
\begin{figure}[htp]
\begin{center}
\begin{tikzpicture}[rotate=90, very thick,scale=1]
\tikzstyle{every node}=[circle, draw=black, fill=white, inner sep=0pt, minimum width=5pt];
    \path (0.1,1.2) node (p1) [label = left: $p_{4}$] {} ;
    \path (2.1,0) node (p4) [label = above right: $p_{1}$]{} ;
    \path (2.3,-1.3) node (p3) [label = right: $p_{2}$] {} ;
     \path (0.3,0) node (p2) [label = below: $p_{3}$] {} ;
        \draw (p1) -- (p4);
      \draw (p3) -- (p4);
     \draw (p2) -- (p3);
      \draw (p2) -- (p1);
            \draw [ultra thick, ->, black!40!white](p1) -- (0.52,0.74);
      \draw [ultra thick, ->, black!40!white](p3) -- (1.88,-0.84);
      \draw [ultra thick, ->, black!40!white](p2) -- (-0.5,-0.65);
      \draw [ultra thick, ->, black!40!white](p4) -- (2.9,0.55);
      \filldraw[fill=black, draw=black]
    (1.2,-0.05) circle (0.004cm);
     \node [draw=white, fill=white] (b) at (-1,-0) {(a)};
              \end{tikzpicture}
    \hspace{0.1cm}
            \begin{tikzpicture}[very thick,scale=1]
\tikzstyle{every node}=[circle, draw=black, fill=white, inner sep=0pt, minimum width=5pt];
    \path (-0.7,0.8) node (p1) [label = left: $p_{1}$] {} ;
    \path (0.7,0.8) node (p4) [label = right: $p_{4}$]{} ;
    \path (-1.6,-0.8) node (p2) [label = below: $p_{2}$] {} ;
     \path (1.6,-0.8) node (p3) [label = below: $p_{3}$] {} ;
      \draw (p1) -- (p4);
    \draw (p1) -- (p2);
    \draw (p3) -- (p4);
    \draw (p2) -- (p3);
     \draw [dashed, thin] (0,-1.4) -- (0,1.4);
     \draw [ultra thick, ->, black!40!white] (p1) -- (-0.7,1.3);
      \draw [ultra thick, ->, black!40!white] (p4) -- (0.7,1.3);
      \draw [ultra thick, ->, black!40!white] (p2) -- (-1.6,-0.3);
      \draw [ultra thick, ->, black!40!white] (p3) -- (1.6,-0.3);
      \node [draw=white, fill=white] (b) at (0,-2.1) {(b)};
        \end{tikzpicture}
        \hspace{0.1cm}
        \begin{tikzpicture}[very thick,scale=1]
\tikzstyle{every node}=[circle, draw=black, fill=white, inner sep=0pt, minimum width=5pt];
    \path (-0.7,0.8) node (p1) [label = left: $p_{1}$]  {} ;
    \path (0.7,0.8) node (p4)[label = right: $p_{4}$] {} ;
    \path (-1.6,-0.8) node [label = below: $p_{2}$](p2)  {} ;
     \path (1.6,-0.8) node [label = below: $p_{3}$](p3)  {} ;
      \draw (p1) -- (p4);
    \draw (p1) -- (p2);
    \draw (p3) -- (p4);
    \draw (p2) -- (p3);
     \draw [dashed, thin] (0,-1.4) -- (0,1.4);
     \draw [ultra thick, ->, black!40!white] (p1) -- (-0.5,0.3);
      \draw [ultra thick, ->, black!40!white] (p4) -- (0.9,1.3);
      \draw [ultra thick, ->, black!40!white] (p2) -- (-2.2,-1);
      \draw [ultra thick, ->, black!40!white] (p3) -- (1,-0.6);
      \node [draw=white, fill=white] (b) at (0,-2.1) {(c)};
        \end{tikzpicture}
\end{center}
\vspace{-0.3cm}
\caption{Infinitesimal motions of frameworks in the plane: (a) a  $\mathcal{C}_2$-symmetric non-trivial infinitesimal motion of $(K_{2,2},p) \in \mathscr{R}_{(K_{2,2},\mathcal{C}_2, \Phi)}$; (b) a $\mathcal{C}_s$-symmetric trivial infinitesimal motion of $(K_{2,2},p) \in \mathscr{R}_{(K_{2,2},\mathcal{C}_s, \Phi_b)}$; (c) a non-trivial infinitesimal motion of $(K_{2,2},p) \in \mathscr{R}_{(K_{2,2},\mathcal{C}_s, \Phi_b)}$ which is not $\mathcal{C}_s$-symmetric.}
\label{fulsym}
\end{figure}

An infinitesimal motion $u$ of a framework $(G,p)\in \mathscr{R}_{(G,S,\Phi)}$  is called \emph{$S$-symmetric} if \begin{equation}\label{fulsymmot} x(u_i)=u_{x(i)}\textrm{ for all } i\in V(G)\textrm{ and all } x\in S\textrm{,}\end{equation} i.e., if $u$ is unchanged under all symmetry operations in $S$ (see also Figure \ref{fulsym}(a) and (b)).

Note that it follows immediately from (\ref{fulsymmot}) that if $u$ is an $S$-symmetric infinitesimal motion of $(G,p)$, then $u_i$ is an element of $U(p_i)$ for each $i$. Moreover, $u$ is uniquely determined by the velocity vectors $u_1,\ldots, u_k$ whenever  $\mathscr{O}_{V(G)}=\{1,\ldots, k\}$ is a set of representatives for the vertex orbits  $S(i)=\{\Phi(x)(i)|\,x\in S\}$ of $G$.

\begin{figure}[htp]
\begin{center}
 \begin{tikzpicture}[very thick,scale=1]
\tikzstyle{every node}=[circle, draw=black, fill=white, inner sep=0pt, minimum width=5pt];
           \node [draw=white, fill=white] (a) at (0.2,-1.0) {$\beta$};
     \node [draw=white, fill=white] (a) at (1,0.2) {$\beta$};
     \node [draw=white, fill=white] (a) at (-1,0.2) {$\beta$};
      \node [draw=white, fill=white] (a) at (-0.1,-0.4) {$\gamma$};
     \node [draw=white, fill=white] (a) at (0.4,0.2) {$\gamma$};
     \node [draw=white, fill=white] (a) at (-0.4,0.2) {$\gamma$};
        \node [draw=white, fill=white] (a) at (0.8,-0.3) {$\alpha$};
     \node [draw=white, fill=white] (a) at (-0.8,-0.3) {$\alpha$};
     \node [draw=white, fill=white] (a) at (-0.2,0.8) {$\alpha$};
        \path (90:1.5cm) node (p6) [label = left: $p_1$] {} ;
    \path (210:1.5cm) node (p7) [label = left: $p_2$] {} ;
    \path (330:1.5cm) node (p8) [label = right: $p_3$] {} ;
     \path (90:0.6cm) node (p9) [label = right: $p_4$] {} ;
     \path (210:0.4cm) node (p10) [label = below: $p_5$] {} ;
     \path (330:0.4cm) node (p11) [label = below: $p_6$] {} ;
   \draw (p6) -- (p9);
      \draw (p8) -- (p11);
     \draw (p7) -- (p8);
      \draw (p7) -- (p6);
     \draw (p6) -- (p8);
      \draw (p7) -- (p10);
            \draw (p9) -- (p10);
      \draw (p11) -- (p10);
          \draw (p9) -- (p11);
          \draw [dashed, thin] (0,-1.4) -- (0,1.9);
           \draw [dashed, thin] (210:2cm) -- (30:2cm);
           \draw [dashed, thin] (330:2cm) -- (150:2cm);
      \node [draw=white, fill=white] (b) at (0,-2.1) {(a)};
         \end{tikzpicture}
        \hspace{2cm}
        \begin{tikzpicture}[very thick,scale=1]
\tikzstyle{every node}=[circle, draw=black, fill=white, inner sep=0pt, minimum width=5pt];
    \path (0,1.4) node (p1) [label = above right: $p_1$] {} ;
    \path (-1,0.4) node (p2) [label = left: $p_2$] {} ;
    \path (1,0.4) node (p3) [label = right: $p_3$] {} ;
   \path (-1,-0.6) node (p4) [label = left: $p_4$]{} ;
   \path (1,-0.6) node (p5) [label = right: $p_5$] {} ;
      \draw (p1) -- (p4);
     \draw (p1) -- (p5);
     \draw (p1) -- (p2);
     \draw (p1) -- (p3);
     \draw (p2) -- (p4);
     \draw (p3) -- (p5);
     \draw (p2) -- (p5);
     \draw (p3) -- (p4);
       \draw [dashed, thin] (0,-1.4) -- (0,1.9);
        \node [rectangle, draw=white, fill=white] (a) at (-0.7,1) {$\alpha$};
     \node [rectangle, draw=white, fill=white] (a) at (0.9,1) {$-\alpha$};
     \node [rectangle, draw=white, fill=white] (a) at (-0.25,0.45) {$\beta$};
     \node [rectangle, draw=white, fill=white] (a) at (0.2,0.45) {$-\beta$};
     \node [rectangle, draw=white, fill=white] (a) at (0.5,-0.7) {$-\gamma$};
     \node [rectangle,draw=white, fill=white] (a) at (-0.5,-0.7) {$\gamma$};
       \node [rectangle, draw=white, fill=white] (a) at (-1.3,-0.1) {$\delta$};
     \node [rectangle, draw=white, fill=white] (a) at (1.3,-0.1) {$-\delta$};
       \node [draw=white, fill=white] (b) at (0,-2.1) {(b)};
        \end{tikzpicture}
        \end{center}
\vspace{-0.3cm}
\caption{Self-stressed frameworks in the plane: (a) a $\mathcal{C}_{3v}$-symmetric self-stress of $(G,p) \in \mathscr{R}_{(G,\mathcal{C}_{3v}, \Phi)}$; (b) a  self-stress of $(H,p) \in \mathscr{R}_{(H,\mathcal{C}_s, \Psi)}$ which is not $\mathcal{C}_s$-symmetric.}
\label{fulsymstr}
\end{figure}
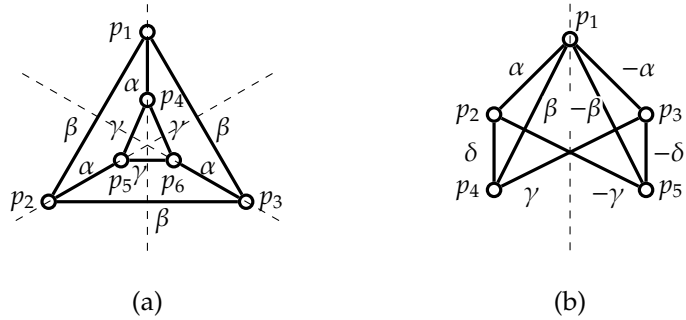

A self-stress $\omega$ of a  framework $(G,p)\in \mathscr{R}_{(G,S,\Phi)}$  is \emph{$S$-symmetric} if $\omega_e=\omega_f$ whenever $e$ and $f$ belong to the same orbit $S(e)=\{\Phi(x)(e)|\,x\in S\}$ of edges of $G$ (see also Figure \ref{fulsymstr}(a)).

Note that an $S$-symmetric self-stress is clearly uniquely determined by the components $\omega_{e_1},\ldots,\omega_{e_r}$, whenever $\mathscr{O}_{E(G)}=\{e_1,\ldots, e_r\}$ is a set of representatives for the edge orbits  $S(e)=\{\Phi(x)(e)|\,x\in S\}$ of $G$.

We are now ready to state the symmetric version of the theorem of Asimov and Roth which was proved in \cite{BS6}.

\begin{theorem} [Schulze \cite{BS6}]
\label{thm:flexes}
Let $G$ be a graph, $S$ be a symmetry group in dimension $d$, $\Phi:S\to \textrm{Aut}(G)$ be a homomorphism, and $(G,p)$ be a framework in $\mathscr{R}_{(G,S,\Phi)}$ whose joints span all of $\mathbb{E}^{d}$. If $(G,p)$ is $S$-regular and $(G,p)$ has an $S$-symmetric non-trivial infinitesimal motion, then there also exists a non-trivial finite motion of $(G,p)$ which preserves the symmetry of $(G,p)$ throughout the path.
\end{theorem}

By Theorem \ref{thm:flexes}, we may test an $S$-regular framework $(G,p)$ for flexibility by analyzing its $S$-symmetric infinitesimal rigidity properties, which in turn can be done via the `orbit rigidity matrix' which was introduced in \cite{BSWWorbit}.

\begin{defin} \emph{Let $G$ be a graph, $S$ be a symmetry group in dimension $d$, $\Phi:S\to \textrm{Aut}(G)$ be a homomorphism, and $(G,p)$ be a framework in $\mathscr{R}_{(G,S,\Phi)}$. Further, let $\mathscr{O}_{V(G)}=\{1,\ldots, k\}$ be a set of representatives for the orbits  $S(i)=\{\Phi(x)(i)|\,x\in S\}$ of vertices of $G$. We construct the \emph{orbit rigidity matrix} (or in short, \emph{orbit matrix}) $\mathbf{O}(G,p,S)$ of $(G,p)$ so that it has exactly one row for each orbit $S(e)=\{\Phi(x)(e)|\,x\in S\}$ of edges of $G$ and exactly $c_i:=\textrm{dim }\big(U(p_i)\big)$ columns for each vertex $i\in \mathscr{O}_{V(G)}$.\\\indent
Given an edge orbit $S(e)$ of $G$, there are two possibilities for the corresponding row in $\mathbf{O}(G,p,S)$:
\begin{description}
\item[Case 1:] The two end-vertices of the edge $e$ lie in distinct vertex orbits. Then there exists an edge in $S(e)$ that is of the form $\{i,x(j)\}$ for some $x\in S$, where $i,j\in\mathscr{O}_{V(G)}$. Let a basis $\mathscr{B}_{i}$ for $U(p_i)$ and a basis  $\mathscr{B}_{j}$ for $U(p_j)$ be given and let $\mathbf{M}_{i}$ and  $\mathbf{M}_{j}$ be the matrices whose columns are the coordinate vectors of  $\mathscr{B}_{i}$ and $\mathscr{B}_{j}$ relative to the canonical basis of $\mathbb{E}^d$, respectively.
The row we write in $\mathbf{O}(G,p,S)$ is:
        \begin{displaymath}\renewcommand{\arraystretch}{0.8}
    \bordermatrix{  & &  i &  &  j &  \cr & 0 \ldots 0 & (p_j-x(p_j))\mathbf{M}_{i} & 0  \ldots  0 & (p_j-x^{-1}(p_i))\mathbf{M}_{j} & 0  \ldots  0}\textrm{.}
    \end{displaymath}
\item[Case 2:] The two end-vertices of the edge $e$ lie in the same vertex orbit. Then there exists an edge in $S(e)$ that is of the form $\{i,x(i)\}$ for some $x\in S$, where $i\in\mathscr{O}_{V(G)}$. The row we write in $\mathbf{O}(G,p,S)$ is:
\begin{displaymath}
\bordermatrix{ & &  i & \cr & 0  \ldots  0 & (2p_i-x(p_i)-x^{-1}(p_i))\mathbf{M}_{i} & 0  \ldots  0}\textrm{.}
           \end{displaymath}
In particular, if $x(p_i)= x^{-1}(p_i)$, this row becomes
\begin{displaymath}
\bordermatrix{ & &  i & \cr & 0  \ldots  0 & 2(p_i-x(p_i))\mathbf{M}_{i} & 0  \ldots  0}\textrm{.}
        \end{displaymath}
\end{description}}
\end{defin}

\begin{remark}\label{rankrem}
\emph{Note that the rank of the orbit rigidity matrix $\mathbf{O}(G,p,S)$ is clearly independent of the choice of bases  for the spaces $U(p_i)$ (and their corresponding matrices $\mathbf{M}_i$), $i=1,\ldots, k$.}
\end{remark}

\begin{remark}\label{nothingfixedrem}
\emph{If none of the joints of $(G,p)$ are fixed by any non-trivial symmetry operation in $S$, then the orbit rigidity matrix $\mathbf{O}(G,p,S)$ of $(G,p)$ has  $dk=d|\mathscr{O}_{V(G)}|$ columns, and each of the matrices $\mathbf{M}_{i}$ and $\mathbf{M}_{j}$ may be chosen to be the $d\times d$ identity matrix. In this case, the matrix $\mathbf{O}(G,p,S)$ becomes particularly easy to construct (see also Example \ref{c2quadexamporbitmat}.)}
\end{remark}

We illustrate the definition of the orbit rigidity matrix with two examples.

\begin{examp}\label{c2quadexamporbitmat}
Consider the $2$-dimensional framework $(K_{2,2},p)\in\mathscr{R}_{(K_{2,2},\mathcal{C}_2)}$ depicted in Figure \ref{fulsym} (a).
If we denote $p_1=(0,a)$, $p_2=(b,c)$, then $p_3=(0,-a)$, and $p_4=(-b,-c)$, and the rigidity matrix of $(K_{2,2},p)$ is the matrix
\begin{displaymath}\bordermatrix{
                &1&2& {3}=C_2(1) &{4}=C_2(2) \cr
                \{1,2\}&(-b,a-c) &  (b,c-a)  &0\ 0 & 0 \ 0\cr
                 \{1,4\}& (b,a+c) & 0 \ 0  & 0 \ 0 & (-b,-a-c)\cr
             C_2\{1,2\}& 0 \ 0  & 0 \ 0 & (b,c-a) &  (-b,a-c) \cr
             C_2\{1,4\}& 0 \ 0  &  (b,a+c)       & (-b,-a-c) & 0 \ 0}\textrm{.}
\end{displaymath}
The orbit matrix $\mathbf{O}(K_{2,2},p,\mathcal{C}_2)$ of $(K_{2,2},p)$ is the matrix
\begin{displaymath}\bordermatrix{
                &1&2 \cr
                \{1,2\}&(p_1-p_2) &  (p_2-p_1)\cr
                 \{1,C_2(2)\}& \big(p_1-C_2(p_2)\big) & \big(p_2-C_2^{-1}(p_1)\big) \cr
             }=\bordermatrix{
                &1&2 \cr
                &(-b,a-c) &  (b,c-a)\cr
                 & (b,a+c) & (b,c+a) \cr
             }
\end{displaymath}
Clearly, $(K_{2,2},p)$ has a one-dimensional space of $\mathcal{C}_{2}$-symmetric infinitesimal flexes spanned by $u=(-1, 0 , x, y ,1 , 0 , -x ,  -y)^T$, where $x=-\frac{c}{a}$ and $y=\frac{b}{a}$. Note that the vector $(-1, 0 , x, y )^T$ lies in the kernel of the orbit matrix $\mathbf{O}(K_{2,2},p,\mathcal{C}_2)$. $\square$
\end{examp}

\begin{examp}\label{csexamporbitmat}
Similarly,  the orbit rigidity matrix of the $2$-dimensional framework $(G,p)\in\mathscr{R}_{(G,\mathcal{C}_{3v})}$ depicted in Figure \ref{fulsymstr} (a)
 is the matrix
\begin{displaymath}\bordermatrix{
                &1&4 \cr
                \{1,4\}&(p_1-p_4)\mathbf{M}_1 &  (p_4-p_1)\mathbf{M}_4\cr
                 \{1,C_3(1)\}& (2p_1-C_3(p_1)-C_3^2(p_1))\mathbf{M}_1 & 0 \cr
                 \{4,C_3(4)\}& 0 & (2p_4-C_3(p_4)-C_3^2(p_4))\mathbf{M}_4 \cr
             },
\end{displaymath}
where $\mathbf{M}_1=\mathbf{M}_4=\bordermatrix{
                & \cr
                &0\cr
                 & 1 \cr}.$
If we denote  $p_1=(0,a)$ and $p_4=(0,b)$, then this matrix is equal to
\begin{displaymath}
\bordermatrix{
                &1&4 \cr
                &(a-b) &  (b-a)\cr
                 & (3a) & 0 \cr
                 & 0 & (3b) \cr
             }.
 \end{displaymath}
Clearly, $(G,p)$ has a one-dimensional space of $\mathcal{C}_{3v}$-symmetric self-stresses, spanned by $(\alpha, \alpha, \alpha, \beta, \beta, \beta, \gamma, \gamma, \gamma)$, where $\alpha=1$, $\beta=\frac{b-a}{3a}$, and $\gamma=\frac{a-b}{3b}$. Note that for $\omega=(\alpha,\beta,\gamma)$, we have $\omega\mathbf{O}(G,p,\mathcal{C}_{3v})=0$. $\square$
\end{examp}

The key result for the orbit rigidity matrix is the following:

\begin{theorem} [Schulze, Whiteley  \cite{BSWWorbit}] \label{orbitmatrixthm}
Let $(G,p)$ be a framework in $\mathscr{R}_{(G,S,\Phi)}$. Then the solutions to $\mathbf{O}(G,p,S)u = 0$ are isomorphic to the space of $S$-symmetric infinitesimal motions of $(G,p)$.
Moreover, the solutions to $\omega\mathbf{O}(G,p,S) = 0$ are isomorphic to the space of  $S$-symmetric self-stresses of  $(G,p)$.
\end{theorem}

As a consequence of Theorem \ref{orbitmatrixthm} we conclude that a framework $(G,p)\in\mathscr{R}_{(G,S,\Phi)}$ is $S$-regular if and only if
$\textrm{\rm rank\,}\big(\mathbf{O}(G,p,S)\big)={\max }\big\{\textrm{\rm rank\,}\big(\mathbf{O}(G,q,S)\big)|\, q\in \mathscr{P}_{(G,S,\Phi)}\big\}$ (see \cite{BS6} and \cite{BSWWorbit} for details).

\section{Euclidean coning of symmetric frameworks}
\label{sec:euclidcone}

In this section we first present the basic process of coning  a framework to the origin without added symmetry (\S3.1).   The basic result of \S3.1 that coning preserves infinitesimal rigidity and independence is also presented in \cite{WWcones}.  Here we describe an alternative rigidity matrix approach to this result which will permit a direct generalization to symmetric frameworks (\S3.2).  We then move joints back and forth along their rays to the origin (the cone joint) in \S3.3.  Finally in
\S3.4 we apply this process to move the framework onto the sphere.

\begin{figure}
    \begin{center}
 \subfloat[]{\includegraphics[width=1.2in]{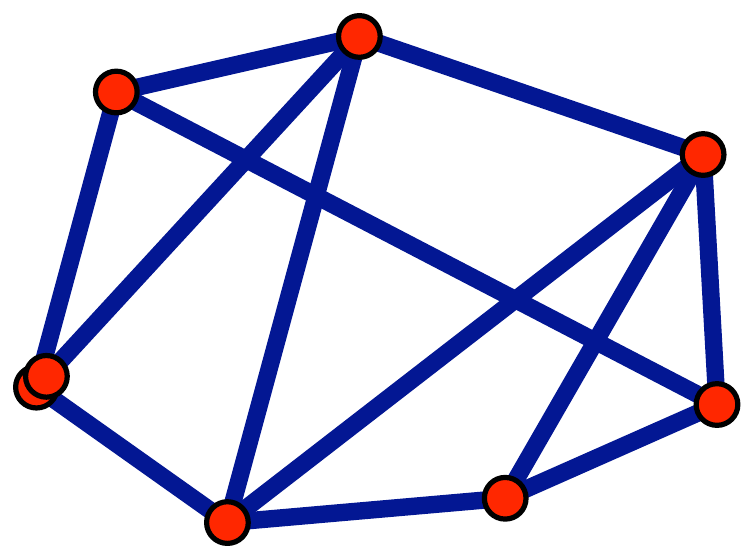}}  \quad
  \subfloat[]{\label{fig:conesimple}\includegraphics[width=1.9in]{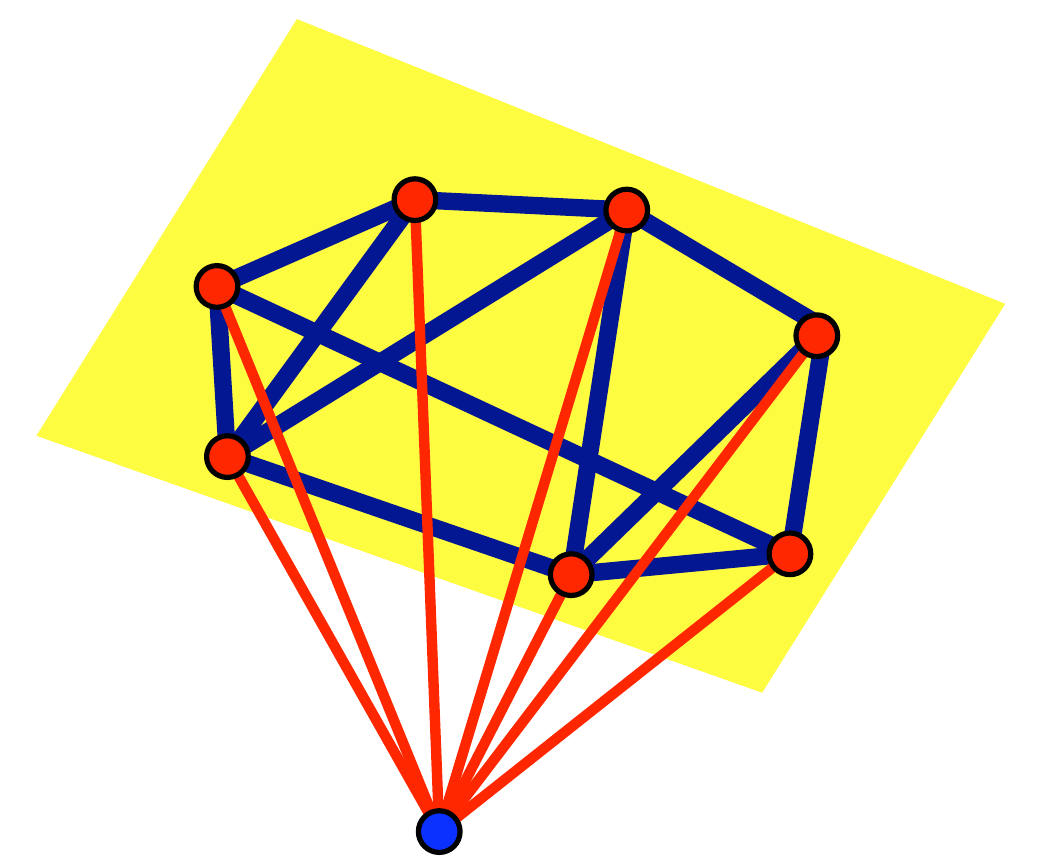}}   \quad
 \subfloat[]{\label{fig:connepulled}\includegraphics[width=1.3in]{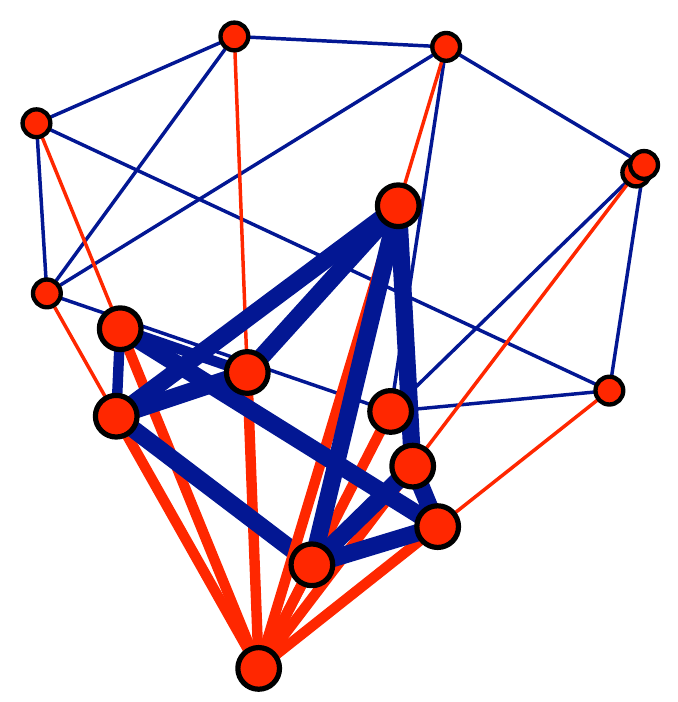}}
      \end{center}
   \caption{The simple cone (b) of a plane framework (a)  has the same static and instantaneous kinematic properties as the more general cone with the same cone rays (c). }
    \label{fig:Cone}
    \end{figure}

In all of this, the matrix processes (row and column reductions) will be reversible, and will permit a simple tracking of how the infinitesimal motions (kernel) and the self-stresses (row dependencies) are modified into isomorphic spaces on the cone.

\subsection{Coning a general framework in $\E^{d}$.}\label{sec:conewithoutsym}

Let $(G,p)$ be a framework in $\E^{d}$.
Recall that the rigidity matrix of $(G,p)$ is the $|E(G)| \times dn$ matrix
\begin{displaymath} \mathbf{R}(G,p)=\bordermatrix{& & & & i & & & & j & & & \cr & & & &  & & \vdots & &  & & &
\cr \{i,j\} & 0 & \ldots &  0 & (p_{i}-p_{j}) & 0 & \ldots & 0 & (p_{j}-p_{i}) &  0 &  \ldots&  0 \cr & & & &  & & \vdots & &  & & &
}
\textrm{.}\end{displaymath}

We embed the framework $(G,p)$ 
into the hyperplane $x_{d+1}=1$  of $\E^{d+1}$ via $\op_{i}=(p_{i},1) \in \E^{d+1}$.  We then cone the resulting framework $(G,\op)$ to the origin $o$ in $\E^{d+1}$, i.e., we connect each of the $n$ joints of $(G,\op)$ with $o$ (see also Figure \ref{fig:Cone}). This adds $n$ rows to the new rigidity matrix  for the \emph{cone framework} $(G*o,\ps)$ in $\E^{d+1}$. We call the underlying graph  $G*o$ of $(G*o,\ps)$ the \emph{cone graph} of $G$.

The
\emph{cone rigidity matrix} of the cone framework $(G*o,\ps)$ is the $(|E(G)|+n) \times (d+1)n$ matrix
\begin{displaymath}{\RC}(G*o,\ps)=\bordermatrix{& & & & i & & & & j & & & \cr & & & &  & & \vdots & &  & & &
\cr
 \{i,j\} & 0 & \ldots &  0 & (\bar{p}_{i}-\bar{p}_{j}) & 0 & \ldots & 0 & (\bar{p}_{j}-\bar{p}_{i}) &  0 &  \ldots&  0
 \cr & & & &  & & \vdots & &  & & &
\cr \{0,i\} & 0 & \ldots &  0 &\bar{p}_{i}& 0 & \ldots & 0 & 0&  0 &  \ldots&  0
\cr & & & &  & & \vdots & &  & & &
\cr  \{0,j\} & 0 & \ldots &  0 & 0& 0 & \ldots & 0 & \bar{p}_{j} &  0 &  \ldots&  0
\cr & & & &  & & \vdots & &  & & &
}
\textrm{.}\end{displaymath}
We have added $n$ rows and $n$ columns. We note that $(\bar{p}_{i}-\bar{p}_{j}) $ is zero in each added column.  Moreover, for each added column (under vertex $i$) there is exactly one added row which is non-zero in this column: $\bar{p}_{i}$ has a 1 in this column.  Thus we have increased the rank by $n$, and preserved the dimension of the kernel.

We have not added the columns for the added cone vertex $o$ - which would have increased the dimension of the kernel by $d+1$. Instead, for convenience, we have introduced a modified cone rigidity matrix, where the kernel is reduced and the infinitesimal motions are restricted to those which fix the origin - the cone joint.  This will be particularly convenient for the sphere (see \S\ref{sec:hemisphere}).

\begin{examp}\label{c2quadexampconeorbit}
Consider the cone framework $(K_{2,2}*o,\ps)\in\mathscr{R}_{(K_{2,2}*o,\mathcal{C}_2^*)}$ of the framework $(K_{2,2},p)\in\mathscr{R}_{(K_{2,2},\mathcal{C}_2)}$ from Example \ref{c2quadexamporbitmat}
(see also Figure \ref{fulsym} (a)). The cone rigidity matrix of $(K_{2,2}*o,\ps)$ is the matrix
\begin{displaymath}\bordermatrix{
                &1&2& 3 &4 \cr
                \{1,2\}&(-b,a-c,0) &  (b,c-a,0)  &0\ 0\ 0 & 0 \ 0\ 0\cr
                 \{1,4\}& (b,a+c,0) & 0 \ 0 \ 0 & 0 \ 0\ 0 & (-b,-a-c,0)\cr
             \{3,4\}&0\ 0\ 0  & 0\ 0\ 0& (b,c-a,0) &  (-b,a-c,0) \cr
             \{2,3\}& 0\ 0\ 0 &  (b,a+c,0)       & (-b,-a-c,0) & 0\ 0\ 0\cr
             \{0,1\}& (0,a,1) &   0\ 0\ 0    & 0\ 0\ 0 & 0\ 0\ 0\cr
             \{0,2\}&  0\ 0\ 0 &   (b,c,1)    & 0\ 0\ 0 & 0\ 0\ 0\cr
             \{0,3\}& 0\ 0\ 0 &   0\ 0\ 0    & (0,-a,1) & 0\ 0\ 0\cr
             \{0,4\}&  0\ 0\ 0 &   0\ 0\ 0    & 0\ 0\ 0 & (-b,-c,1)\cr
             }.
\end{displaymath}$\square$
\end{examp}

If we apply this same process to the complete graph on the vertices, we see that the trivial motions in $\E^{d}$ go to trivial motions in $\E^{d+1}$ which fix the origin.  Rotations will go to rotations around extended axes which contain the origin, and translations will also go to rotations around axes formed by joining the origin to the     `implied center of the translation' at infinity.

\begin{figure}
    \begin{center}
 \subfloat[]{\includegraphics[width=1.4in]{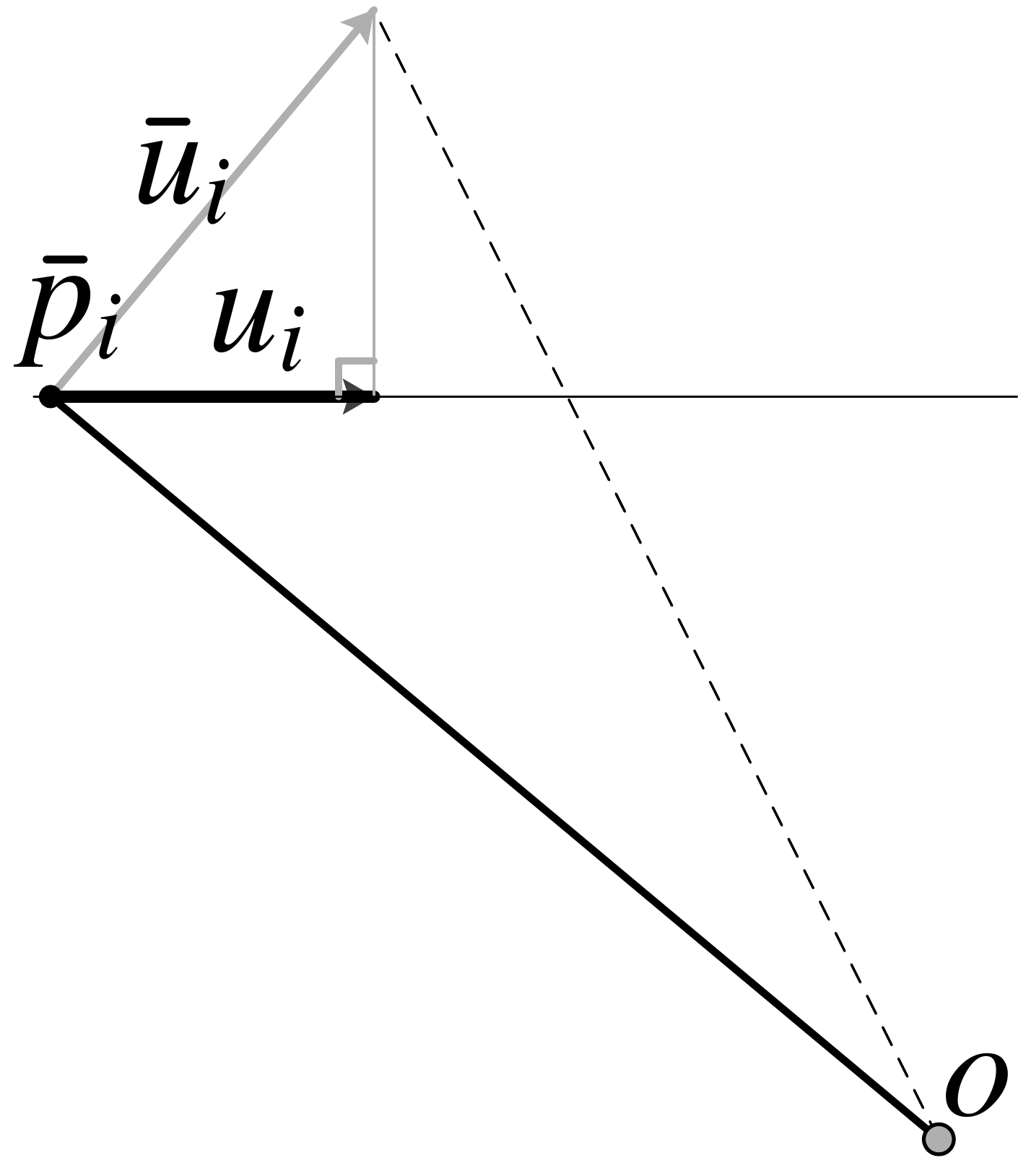}}  \quad
  \subfloat[]{\label{fig:conevelocity}\includegraphics[width=1.4in]{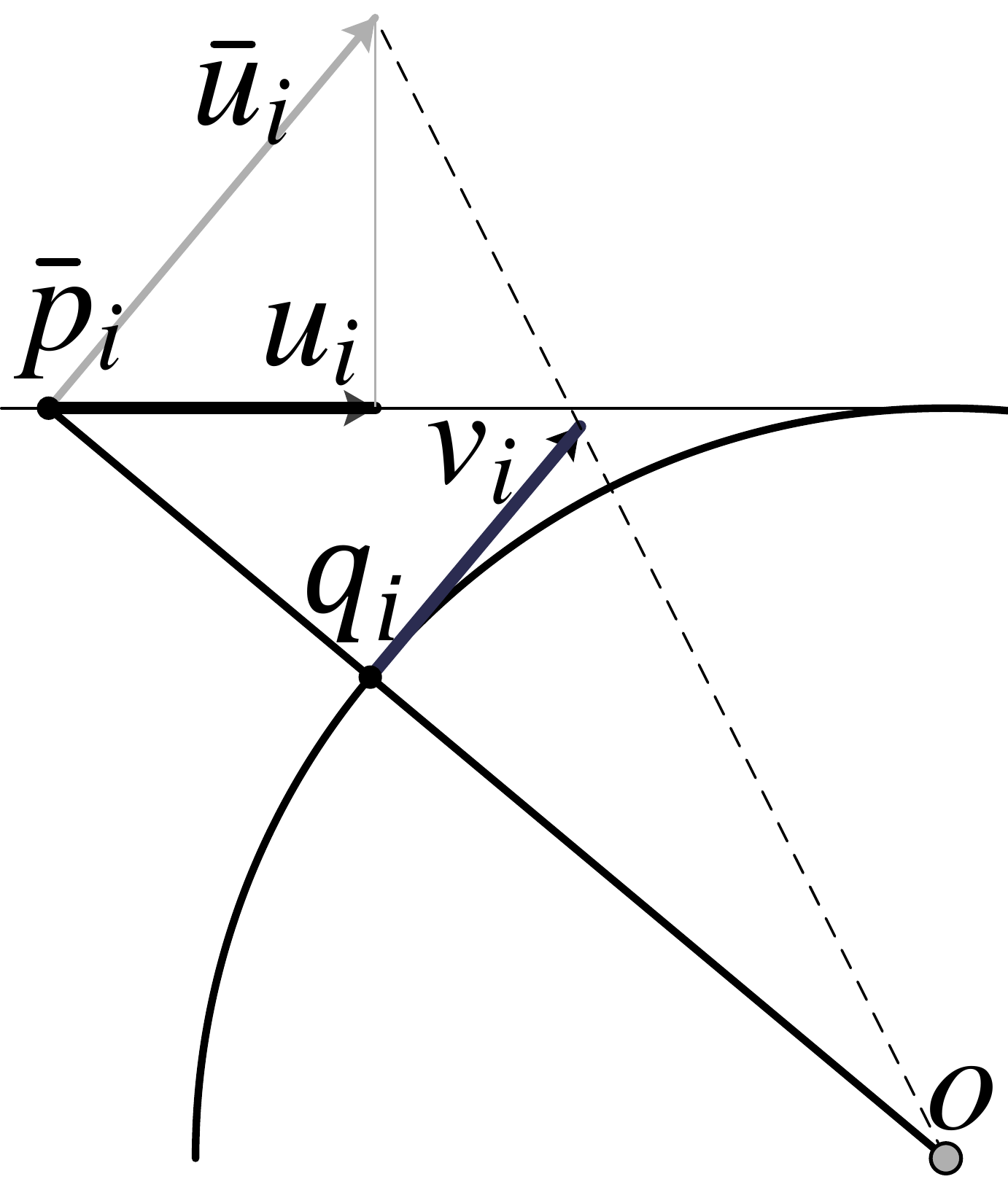}}   \quad
 \subfloat[]{\label{fig:spherevelocity}\includegraphics[width=1.5in]{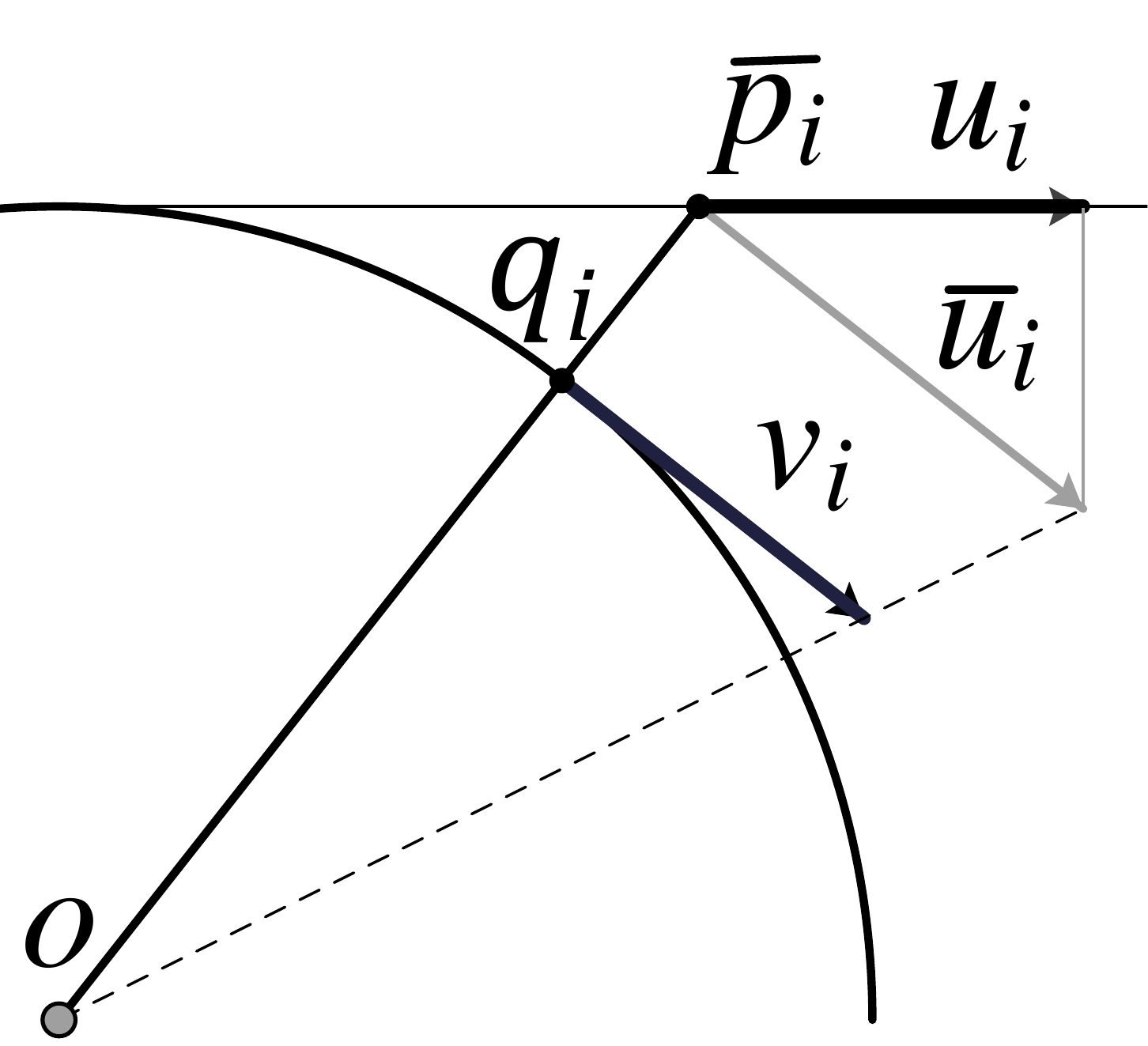}}
      \end{center}
  \caption{The velocity in Euclidean space $u_{i}$ goes to a distinct velocity $\bar{u}_{i}$ on the cone (a).    This velocity is modified when pulled onto the hemisphere (b), (c) but if $u_i=u_j$ then $||v_i|=||v_j||$. }
    \label{fig:SpherePlane}
    \end{figure}

What does the transfer of velocities look like, in practice?  Consider Figure~\ref{fig:SpherePlane} which shows the process from the line to the plane.
We keep the same first $d$ coordinates of the velocities, and add whatever last entry will make the vector perpendicular to the bar from the origin to the joint.  That is, we take the unique vector which is perpendicular to the coning bar and projects orthogonally onto the previous velocity. Explicitly, for a $d$-dimensional velocity vector $u_i$ at a joint $p_i$, the new velocity vector at the joint $\bar p_i$ is the $(d+1)$-dimensional vector $(u_i,-u_i\cdot p_i)$.

If we consider self-stresses - row dependencies - the new rows are independent of all previous rows, and there is a direct correspondence of self-stresses between the two matrices (simply add $0$ coefficients to all the new rows when extending, or delete these coefficients when transferring back).

\begin{theorem} [General coning] \label{thm:ConingGeneral}
Let $(G,p)$ be a framework in $\mathbb{E}^d$.
\begin{enumerate}
\item[(i)] The space of infinitesimal motions of $(G,p)$ in $\E^{d}$ is isomorphic to the space of infinitesimal motions of the cone framework
$(G*o,\ps)$ in $\E^{d+1}$, with the cone joint  fixed at the origin;
 \item[(ii)]  the space of trivial infinitesimal motions of $(G,p)$ in $\E^{d}$ is isomorphic to the space of trivial infinitesimal motions of
$(G*o,\ps)$ in $\E^{d+1}$, with the cone joint  fixed at the origin;
 \item[(iii)]  $(G,p)$ has a non-trivial infinitesimal motion in $\E^{d}$ if and only if $(G*o,\ps)$ has a non-trivial infinitesimal motion in $\E^{d+1}$;
\item[(iv)]  the space of self-stresses of the framework $(G,p)$  is isomorphic to the space of self-stresses of $(G*o,\ps)$;
\item[(v)]  $(G,p)$ is infinitesimally rigid (isostatic) in $\E^{d}$ if and only if $(G*o,\ps)$ is infinitesimally rigid (isostatic) in $\E^{d+1}$.
\end{enumerate}
\end{theorem}

Notice that we have stated the last three parts above without any reference to the cone joint being fixed.   This was intentional.   Adding the missing $d+1$ columns for the cone vertex will add a $(d+1)$-dimensional space of infinitesimal translations, but will not add any non-trivial infinitesimal motions.  The rank of the matrix will not change, nor will the space of self-stresses.  It is a general principle of statics that an equilibrium at all but one vertex will satisfy the equilibrium equation at the final vertex.
\medskip

As a next step, for each $i$, we subtract the added row  for the coning edge $\{0,i\}$ from each edge containing the vertex $i$. This yields the following modified cone rigidity matrix:
\begin{displaymath} {\Rb}(G*o,\ps)=\bordermatrix{& & & & i & & & & j & & &
\cr & & & &  & & \vdots & &  & & &
\cr
 \{i,j\} & 0 & \ldots &  0 & (-\bar{p}_{j}) & 0 & \ldots & 0 & (-\bar{p}_{i}) &  0 &  \ldots&  0
 \cr & & & &  & & \vdots & &  & & &
\cr \{0,i\} & 0 & \ldots &  0 &\bar{p}_{i}& 0 & \ldots & 0 & 0&  0 &  \ldots&  0
\cr & & & &  & & \vdots & &  & & &
\cr  \{0,j\} & 0 & \ldots &  0 & 0& 0 & \ldots & 0 & \bar{p}_{j} &  0 &  \ldots&  0
\cr & & & &  & & \vdots & &  & & &
}
\textrm{.}\end{displaymath}
This modified matrix no longer has the direct appearance of a 'rigidity matrix for a framework' - but the form is well structured for further work throughout this section.

This row reduction preserves the rank and the solution space, so that the kernel (the space of infinitesimal motions) remains unchanged.  The row dependencies (self-stresses), however, do take a different form.

Specifically, the former equilibrium equation for each vertex $i$:   $\Sigma_{j\ | \ \{i,j\}\in E(G)}\omega_{ij} (\op_{i} - \op_{j})$ is rewritten in the form $\Sigma_{j\ | \ ij\in E(G)}\omega_{ij} (- \op_{j})  + (\Sigma_{j\ | \ \{i,j\}\in E(G)}\omega_{ij} )\op_{i}$.  So we keep the same coefficient on the rows for the edges $\{i,j\}$ and  the new coefficient $(\Sigma_{j\ | \ \{i,j\}\in E(G)}\omega_{ij} )$ on the row for $\{0,i\}$.  We can now see the direct isomorphism of the spaces of dependencies between the two matrices.

Since we no longer have the explicit rigidity matrix for an identified framework, we record these correspondences in the following form.

\begin{prop} For a given cone framework $(G*o,\op)$ the matrices $\bar{R}(G*o,\op)$ and $\RC(G*o,\op)$ have have the same rank as well as isomorphic  kernels and cokernels.
\end{prop}

\subsection{Coning with symmetry}\label{sec:conewithsym}

In this section, we show that for a given framework $(G,p)$ in $\mathscr{R}_{(G,S,\Phi)}$, the cone framework $(G*o,\ps)$ has  essentially the same symmetry group and the same symmetric infinitesimal rigidity properties as $(G,p)$. Since the type $\Phi:S\to \textrm{Aut}(G)$ of a given framework $(G,p)$ with symmetry group $S$ is always uniquely determined (recall Section \ref{sec:sym}), we will simply write $(G,p)\in\mathscr{R}_{(G,S)}$ from now on.

\begin{defin} \emph{Let $S$ be a symmetry group in dimension $d$ and let $x\in S$. We denote the matrix which represents $x$ with respect to the canonical basis of $\mathbb{E}^d$ by $M_x$. Then $x^*$ is defined as the orthogonal transformation of $\mathbb{E}^{d+1}$ which is represented by the matrix $$M_{x^*}=\left( \begin{array}{cccc} & & & 0\\& M_x & & \vdots\\ & & & 0\\0 & \cdots & 0 &1 \end{array} \right)$$
with respect to the canonical basis of $\mathbb{E}^{d+1}$. Further, we let $S^*$ be the orthogonal group in dimension $(d+1)$ which has the elements $\{x^*:x\in S\}$.
}
\end{defin}

\begin{theorem}
\label{thm:groupextend}
A framework $(G,p)$ is an element of $\mathscr{R}_{(G,S)}$ if and only if the cone framework $(G*o,\ps)$ is an element of $\mathscr{R}_{(G*o,S^*)}$.
\end{theorem}
\textbf{Proof.} Let $x\in S$. We have $M_x p_i^T=p^T_{x(i)}$ for all $i\in V(G)$ and all $x\in S$ if and only if
$$M_{x^*} \bar{p}_i^T=\left( \begin{array}{cccc} & & & 0\\& M_x & & \vdots\\ & & & 0\\0 & \cdots & 0 &1 \end{array} \right)\left( \begin{array}{c} \\ p_i^T \\ \\1 \end{array} \right)= \left( \begin{array}{c} \\ M_x p_i^T \\ \\1 \end{array} \right)= \left( \begin{array}{c} \\ p_{x(i)}^T \\ \\1 \end{array} \right)=\bar{p}_{x^*(i)}^T$$
for all $i\in V(G)$ and all $x\in S$. This gives the result. $\square$

\begin{theorem} \label{th:rigpropcone} Let $(G,p)$ be a framework in $\mathscr{R}_{(G,S)}$, and $(G*o,\ps)$ be the corresponding cone framework in
$\mathscr{R}_{(G*o,S^*)}$. Then:
\begin{enumerate}
\item[(i)] the space of $S$-symmetric infinitesimal motions of $(G,p)$ in $\E^{d}$ is isomorphic to the space of $S^*$-symmetric infinitesimal motions of $(G*o,\ps)$ in $\E^{d+1}$, with the cone joint of $(G*o,\ps)$ fixed at the origin;
    \item[(ii)]  the space of trivial $S$-symmetric infinitesimal motions of $(G,p)$ in $\E^{d}$ is isomorphic to the space of trivial $S^*$-symmetric infinitesimal motions of $(G*o,\ps)$ in $\E^{d+1}$, with the cone joint of $(G*o,\ps)$ fixed at the origin;
    \item[(iii)]   $(G,p)$ has a non-trivial $S$-symmetric infinitesimal motion in $\mathbb{E}^d$ if and only if $(G*o,\ps)$ has a non-trivial $S^*$-symmetric infinitesimal motion in $\mathbb{E}^{d+1}$;
\item[(iv)]  the space of $S$-symmetric self-stresses of $(G,p)$ is isomorphic to the space of $S^*$-symmetric self-stresses of $(G*o,\ps)$.
 \end{enumerate}
\end{theorem}
\textbf{Proof.} We use the same basic procedure as in Section \ref{sec:conewithoutsym}.
Recall that for a $d$-dimensional framework $(G,p)\in\mathscr{R}_{(G,S)}$ whose vertex orbits under the action of $S$ are represented by the vertices $1\ldots, k$, the orbit rigidity matrix  $\mathbf{O}(G,p,S)$ of $(G,p)$ is the $r \times c$ matrix
\begin{displaymath}\renewcommand{\arraystretch}{0.7}
    \bordermatrix{  & &  i &  &  j & \cr &&& \vdots &&\cr \{i,x(j)\} & 0 \ldots 0 & (p_i-x(p_j))\mathbf{M}_{i} & 0  \ldots  0 & (p_j-x^{-1}(p_i))\mathbf{M}_{j} & 0  \ldots  0\cr
    &&& \vdots &&\cr
    \{i,y(i)\} & 0  \ldots  0 & (2p_i-y(p_i)-y^{-1}(p_i))\mathbf{M}_{i} & 0  \ldots  0 &0 & 0  \ldots  0\cr
      &&& \vdots &&
    }\textrm{,}
    \end{displaymath}
    where $c=\sum_{i=1}^{k} c_i$.

We embed the framework $(G,p)$ into the hyperplane $x_{d+1}=1$ of $\mathbb{E}^{d+1}$ via $\op_{i}=(p_{i},1) \in \E^{d+1}$.  We then cone this to the origin $o$ in $\E^{d+1}$. This gives rise to $k$ new rows in the orbit rigidity matrix $\mathbf{\OC}(G*o,\ps,S^*)$ of the framework $(G*o,\ps)\in\mathscr{R}_{(G*o,S^*)}$.

For each $i=1,\ldots, k$, we let $\mathscr{B}^*_{i}$ be the basis of the space $U(\ps_i)$ which consists of the basis vectors $\{(w,0): w\in \mathscr{B}_{i}\}$ and the additional basis vector $e_{d+1}$ - the $(d+1)$st canonical basis vector of $\mathbb{E}^{d+1}$. Further, we let $\mathbf{M}_i^*$  be the matrix whose columns are the coordinate vectors of  $\mathscr{B}^*_{i}$ relative to the canonical basis of $\mathbb{E}^{d+1}$; that is, $$\mathbf{M}_i^*=\left( \begin{array}{cccc} & & & 0\\& \mathbf{M}_i & & \vdots\\ & & & 0\\0 & \cdots & 0 &1 \end{array} \right).$$
The orbit rigidity matrix $\mathbf{\OC}(G*o,\ps,S^*)$  has the form
   \begin{displaymath}\renewcommand{\arraystretch}{0.7}
    \bordermatrix{  & &  i &  &  j & \cr &&& \vdots &&\cr \{i,x(j)\} & 0 \ldots 0 & (\bar p_i-x^*(\bar p_j))\mathbf{M}_i^* & 0  \ldots  0 & (\bar p_j-x^{*-1}(\bar p_i))\mathbf{M}_j^* & 0  \ldots  0\cr
    &&& \vdots &&\cr
     \{i,y(i)\}& 0  \ldots  0 & (2 \bar p_i-y^*(\bar p_i)-y^{*-1}(\bar p_i))\mathbf{M}_i^* & 0  \ldots  0 &0 & 0  \ldots  0\cr
      &&& \vdots &&
      \cr \{0,i\}& 0 \ldots  0 &\bar{p}_{i}\mathbf{M}_i^*& 0 \ldots 0 & 0&  0  \ldots  0
\cr & & &  \vdots  & &
\cr \{0,j\}  & 0 \ldots  0 & 0& 0  \ldots  0 & \bar{p}_{j}\mathbf{M}_j^* &  0  \ldots  0
\cr & & &  \vdots  & &
      }\textrm{.}
    \end{displaymath}

Note that we have added $k$ rows and $k$ columns, and that for every edge $\{i,x(j)\}$, $i\neq j$, and every edge $\{i,y(i)\}$ of $G$, the vectors $(\bar p_i-x^*(\bar p_j))\mathbf{M}_i$ and $(2 \bar p_i-y^*(\bar p_i)-y^{*-1}(\bar p_i))\mathbf{M}_i^*$ are zero in the added column.  Moreover, for each added column (under vertex $i$) there is exactly one added row which is non-zero in this column: $\bar{p}_{i}\mathbf{M}_i^*$ has a $1$ in this column.  Thus we have increased the rank by $k$, and preserved the dimension of the kernel.

If we apply this process to the complete graph on the vertices, we see that the dimension $m$ of the space of trivial $S$-symmetric infinitesimal motions of $(G,p)$ is equal to the dimension $m^*$ of the space of trivial $S^*$-symmetric infinitesimal motions of $(G*o,\ps)$ (with the cone joint fixed at the origin). This gives the result. $\square$

Similar to Theorem \ref{thm:ConingGeneral}, we have stated the last two parts in Theorem \ref{th:rigpropcone} without any reference to the cone joint being fixed.  This is because the missing $c_0:=\textrm{dim}(U(p_0))=\textrm{dim}(\bigcap_{x^*\in S^*}F_{x^*})$ columns for the cone joint will add a $c_0$-dimensional space of  $S^*$-symmetric infinitesimal translations, but will not add any $S^*$-symmetric non-trivial infinitesimal motions.  The rank of the matrix will not change, nor will the space of row dependencies (symmetric self-stresses).

The transfer of symmetric velocities from the framework $(G,p)$ to the cone framework $(G*o,\ps)$ is the same as the one described in Section \ref{sec:conewithoutsym}: we keep the same first $d$ coordinates of the velocities, and add whatever last entry will make the vector perpendicular to the bar from the origin to the joint. Note that this transfers an $S$-symmetric infinitesimal motion of $(G,p)$ to an $S^*$-symmetric infinitesimal motion of $(G*o,\ps)$. Similarly, since the new rows in $\mathbf{\OC}(G*o,\ps,S^*)$ are independent of all previous rows, there is a direct correspondence of row-dependencies (symmetric self-stresses) between the two matrices $\mathbf{O}(G,p,S)$ and $\mathbf{\OC}(G*o,\ps,S^*)$ (simply add $0$ coefficients to all the new rows when extending, or delete these coefficients when transferring back).

Note that in the case where $S$ is the identity group, Theorem \ref{th:rigpropcone} simply restates Theorem \ref{thm:ConingGeneral}.

Analogous to the process described in Section \ref{sec:conewithoutsym}, we now carry out a row reduction in the cone orbit matrix $\mathbf{\OC}(G*o,\ps,S^*)$
(using the added row corresponding to coning edge $\{0,i\}$ to subtract from edges containing the vertex $i$). This produces a modified cone orbit matrix $\mathbf{\Ob}(G*o,\ps,S^*)$ which is equal to
   \begin{displaymath}\renewcommand{\arraystretch}{0.7}
    \bordermatrix{  & & i &  &  j & \cr &&& \vdots &&\cr \{i,x(j)\} & 0 \ldots 0 & (-x^*(\bar p_j))\mathbf{M}_i^* & 0  \ldots  0 & (-x^{*-1}(\bar p_i))\mathbf{M}_j^* & 0  \ldots  0\cr
    &&& \vdots &&\cr
     \{i,y(i)\}& 0  \ldots  0 & (\bar p_i-y^*(\bar p_i)-y^{*-1}(\bar p_i))\mathbf{M}_i^* & 0  \ldots  0 &0 & 0  \ldots  0\cr
      &&& \vdots &&
      \cr \{0,i\}& 0 \ldots  0 &\bar{p}_{i}\mathbf{M}_i^*& 0 \ldots 0 & 0&  0  \ldots  0
\cr & & &  \vdots  & &
\cr \{0,j\}  & 0 \ldots  0 & 0& 0  \ldots  0 & \bar{p}_{j}\mathbf{M}_j^* &  0  \ldots  0
\cr & & &  \vdots  & &
      }\textrm{.}
    \end{displaymath}
Analogous to the matrix $\bar{R}(G*o,\op)$ in Section \ref{sec:conewithoutsym}, the matrix $\mathbf{\Ob}(G*o,\ps,S^*)$ no longer has the direct appearance of an 'orbit rigidity matrix for a framework'. However, we will be using this matrix as a starting point for further work throughout this section.

This row reduction preserves the rank and the solution space, so that the kernel (the space of symmetric infinitesimal motions) remains unchanged.  The row dependencies (symmetric self-stresses) however change in the analogous way as described in Section \ref{sec:conewithoutsym}.
Specifically, we keep the same coefficients, $\omega_{ij}$, on the rows for the original edges (i.e., the edges of the form $\{i,x(j)\}$ and $\{i,y(i)\}$) and  the new coefficient $(\Sigma_{j\ | \ \{i,j\}\in \mathscr{O}_{E(G)}}\omega_{ij} )$ on the row for $\{0,i\}$.

We summarize these observations in the following theorem:


\begin{theorem}
For a given cone framework $(G*o,\op)\in\mathscr{R}_{(G*o,S^*)}$, the matrices $\mathbf{\OC}(G*o,\ps,S^*)$ and $\mathbf{\Ob}(G*o,\ps,S^*)$ have have the same rank as well as isomorphic  kernels and cokernels.
\end{theorem}

We conclude this section with two examples, illustrating Theorem \ref{th:rigpropcone}.

\begin{examp}
Consider the framework $(K_{2,2},p)$ from Examples \ref{c2quadexamporbitmat} and \ref{c2quadexampconeorbit}. The orbit rigidity matrix $\mathbf{\OC}(K_{2,2}*o,\ps,\mathcal{C}_2^*)$ of the cone framework $(K_{2,2}*o,\ps)$ is the matrix
\begin{displaymath}\bordermatrix{
                &1&2 \cr
                \{1,2\}&(\bar p_1-\bar p_2) &  (\bar p_2-\bar p_1)\cr
                 \{1,C_2^*(2)\}& \big(\bar p_1-C_2^*(\bar p_2)\big) & \big(\bar p_2-C_2^{*-1}(\bar p_1)\big) \cr
                 \{0,1\}&\bar p_1 &  0\ 0\ 0\cr
                 \{0,2\}& 0\ 0\ 0 &  \bar p_2\cr
             }=\bordermatrix{
                &1&2 \cr
                &(-b,a-c,0) &  (b,c-a,0)\cr
                 & (b,a+c,0) & (b,c+a,0) \cr
                 & (0,a, 1) &  0\ 0\ 0\cr
                 & 0\ 0\ 0 &  (b,c,1)\cr
             }
\end{displaymath}
and the modified matrix $\mathbf{\Ob}(K_{2,2}*o,\ps,\mathcal{C}_2^*)$ is the matrix
\begin{displaymath}\bordermatrix{
                &1&2 \cr
                \{1,2\}&(-\bar p_2) &  (-\bar p_1)\cr
                 \{1,C_2^*(2)\}& \big(-C_2^*(\bar p_2)\big) & \big(-C_2^{*-1}(\bar p_1)\big) \cr
                 \{0,1\}&\bar p_1 &  0\ 0\ 0\cr
                 \{0,2\}& 0\ 0\ 0 &  \bar p_2\cr
             }=\bordermatrix{
                &1&2 \cr
                &(-b,-c,-1) &  (0,-a,-1)\cr
                 & (b,c,-1) & (0,a,-1) \cr
                 & (0,a, 1) &  0\ 0\ 0\cr
                 & 0\ 0\ 0 &  (b,c,1)\cr
             }.
\end{displaymath}
The $\mathcal{C}_2^*$-symmetric infinitesimal flex of $(K_{2,2}*o,\ps)$ corresponding to the $\mathcal{C}_2$-symmetric infinitesimal flex $u=(-1, 0 , x, y ,1 , 0 , -x ,  -y)^T$ of $(K_{2,2},p)$ defined in Example \ref{c2quadexamporbitmat} is $u^*=(-1, 0 ,0, x, y ,0,1 , 0 , 0, -x ,  -y, 0)^T$. (Note that the vector $(-1, 0 ,0, x, y ,0)^T$ lies in the kernel of the matrix $\mathbf{\OC}(K_{2,2}*o,\ps,\mathcal{C}_2^*)$.) $\square$
\end{examp}


\begin{examp}\label{c3k4examp}
Consider the $2$-dimensional framework $(G,p)\in\mathscr{R}_{(G,\mathcal{C}_{3v})}$ from Example \ref{csexamporbitmat} which is depicted in Figure \ref{fulsymstr} (a). The orbit rigidity matrix $\mathbf{\OC}(G*o,\ps,\mathcal{C}_{3v}^*)$ of the cone framework $(G*o,\ps)$ is the matrix
\begin{displaymath}\bordermatrix{
                &1&4 \cr
                \{1,4\}&(\bar p_1-\bar p_4)\mathbf{M}_1^* &  (\bar p_4-\bar p_1)\mathbf{M}_4^*\cr
                 \{1,C_3^*(1)\}& (2\bar p_1-C_3^*(\bar p_1)-C_3^{2*}(\bar p_1))\mathbf{M}_1^* & 0\ 0 \cr
                 \{4,C_3^*(4)\}& 0\ 0 & (2\bar p_4-C_3^*(\bar p_4)-C_3^{2*}(\bar p_4))\mathbf{M}_4^* \cr
                 \{0,1\}& \bar{p}_{1}\mathbf{M}_1^* & 0\ 0 \cr
                 \{0,4\}& 0\ 0 & \bar{p}_{4}\mathbf{M}_4^* \cr
             }.
\end{displaymath}
With the notation of Example \ref{csexamporbitmat}, this matrix is equal to
\begin{displaymath}
\bordermatrix{
                &1&4 \cr
                &(a-b,0) &  (b-a,0)\cr
                 & (3a,0) & 0\ 0 \cr
                 & 0\ 0 & (3b,0) \cr
                & (a,1) & 0\ 0 \cr
                 & 0\ 0 & (b,1) \cr
             }.
 \end{displaymath}
The $\mathcal{C}_{3v}^*$-symmetric self-stress of $(G*o,\ps)$ corresponding to the $\mathcal{C}_{3v}$-symmetric self-stress $\omega$ of $(G,p)$ defined in Example \ref{csexamporbitmat} is $\omega^*=(\alpha,\beta,\gamma,0,0)$.
Finally, note that the modified matrix $\mathbf{\Ob}(G*o,\ps,\mathcal{C}_{3v}^*)$ is the matrix
\begin{displaymath}
\bordermatrix{
                &1&4 \cr
                &(-b,-1) &  (-a,-1)\cr
                 & (2a,-1) & 0\ 0 \cr
                 & 0\ 0 & (2b,-1) \cr
                 & (a,1) & 0\ 0 \cr
                 & 0\ 0 & (b,1) \cr
             }.
  \end{displaymath}
A row dependency of this matrix is given by the vector $(\alpha,\beta,\gamma,\alpha+\beta,\alpha+\gamma)$. $\square$
\end{examp}


\subsection{Pulling vertex orbits}

If we move an orbit of joints off of the hyperplane $x_{d+1}=1$, along the corresponding rays to the cone joint, then this amounts to multiplying the coordinates $\bar{p}_{i}$ by a scalar $\alpha_{i}\neq 0$ to create $\bar{q}_{i}:=\alpha_{i}\bar{p}_{i}$.  In the orbit matrix $\mathbf{\Ob}(G*o,\ps,S^*)$, we multiply the rows $\{i,x(j)\}$, $i\neq j$, by $\alpha_i\alpha_j$, and the rows $\{i,y(i)\}$ and $\{0,i\}$ by $\alpha_{i}^{2}$ to produce
   \begin{displaymath}\renewcommand{\arraystretch}{0.7}
    \bordermatrix{  & &  i &  &  j & \cr &&& \vdots &&\cr \{i,x(j)\} & 0 \ldots 0 & \alpha_i\alpha_j(-x^*(\bar p_j))\mathbf{M}_i^* & 0  \ldots  0 & \alpha_i\alpha_j(-x^{*-1}(\bar p_i))\mathbf{M}_j^* & 0  \ldots  0\cr
    &&& \vdots &&\cr
     \{i,y(i)\}& 0  \ldots  0 & \alpha_i^2(\bar p_i-y^*(\bar p_i)-y^{*-1}(\bar p_i))\mathbf{M}_i^* & 0  \ldots  0 &0 & 0  \ldots  0\cr
      &&& \vdots &&
      \cr \{0,i\}& 0 \ldots  0 &\alpha_i^2\bar{p}_{i}\mathbf{M}_i^*& 0 \ldots 0 & 0&  0  \ldots  0
\cr & & &  \vdots  & &
\cr \{0,j\}  & 0 \ldots  0 & 0& 0  \ldots  0 & \alpha_j^2\bar{p}_{j}\mathbf{M}_j^* &  0  \ldots  0
\cr & & &  \vdots  & &
      }\textrm{.}
    \end{displaymath}
This is equivalent to multiplying the orbit matrix $\mathbf{\Ob}(G*o,\ps,S^*)$ on the left by an invertible  matrix.
We can then multiply the columns for vertex $i$ by ${1\over \alpha_{i}}$ to get:
   \begin{displaymath}\renewcommand{\arraystretch}{0.7}
    \bordermatrix{  & &  i &  &  j & \cr &&& \vdots &&\cr \{i,x(j)\} & 0 \ldots 0 &(-x^*(\bar q_j))\mathbf{M}_i^* & 0  \ldots  0 & (-x^{*-1}(\bar q_i))\mathbf{M}_j^* & 0  \ldots  0\cr
    &&& \vdots &&\cr
     \{i,y(i)\}& 0  \ldots  0 & (\bar q_i-y^*(\bar q_i)-y^{*-1}(\bar q_i))\mathbf{M}_i^* & 0  \ldots  0 &0 & 0  \ldots  0\cr
      &&& \vdots &&
      \cr \{0,i\}& 0 \ldots  0 &\bar{q}_{i}\mathbf{M}_i^*& 0 \ldots 0 & 0&  0  \ldots  0
\cr & & &  \vdots  & &
\cr \{0,j\}  & 0 \ldots  0 & 0& 0  \ldots  0 & \bar{q}_{j}\mathbf{M}_j^* &  0  \ldots  0
\cr & & &  \vdots  & &
      }\textrm{.}
    \end{displaymath}
This is equivalent to multiplying the orbit matrix on the right by an invertible  matrix.

Since all of these changes are reversible equivalences, representable by invertible matrix multiplication, we see that the general cone framework has no $S^*$-symmetric infinitesimal
flex  ($S^*$-symmetric self-stress) if and only if the projection onto a  hyperplane has no $S$-symmetric infinitesimal flex ($S$-symmetric self-stress).
More generally, coning gives an isomorphism of the spaces of $S^*$-symmetric self-stresses of the coned framework and $S$-symmetric self-stresses of the framework projected into a hyperplane
not containing the cone joint,
 as well as an isomorphism of the space of $S^*$-symmetric first-order motions of the coned framework which fix the cone joint, and the $S$-symmetric first-order motions of the projected framework,
an isomorphism which takes the trivial motions to the trivial motions.

\begin{theorem}[Transfer for coning with symmetry]\label{thmtranconwithsym}
Let $q$ be a configuration of $n+1$ points (including the origin) in $\E^{d+1}$ such that the projection $\pi(q)$ from the origin onto the hyperplane and then projected back to $\E^{d}$, is equal to $p \in \E^{dn}$. Then
\begin{enumerate}
\item[(i)]  the space of $S$-symmetric infinitesimal motions  of $(G,p)$ in $\E^{d}$ is isomorphic to the space of $S^{*}$-symmetric infinitesimal motions of $(G*o,q)$ in $\E^{d+1}$, with the cone joint of $(G*o,q)$ fixed at the origin;
\item[(ii)] the space of $S$-symmetric trivial infinitesimal motions  of $(G,p)$ in $\E^{d}$ is isomorphic to the space of $S^{*}$-symmetric trivial infinitesimal motions of $(G*o,q)$ in $\E^{d+1}$, with the cone joint of $(G*o,q)$ fixed at the origin;
\item[(iii)]   $(G,p)$ has a non-trivial $S$-symmetric  infinitesimal motion in $\E^{d}$  if and only if $(G*o,q)$ has a non-trivial $S^{*}$-symmetric  infinitesimal motion in $\E^{d+1}$;
\item[(iv)]  the space of $S$-symmetric  self-stresses  of $(G,p)$ is isomorphic to the space of $S^{*}$-symmetric  self-stresses of $(G*o,q)$.
\end{enumerate}
\end{theorem}

For $S=Id$, this is the original result given in \cite{WWcones}, using different techniques.

\subsection{Hemispherical realizations}\label{sec:hemisphere}
In particular, we could scale all joints onto the sphere of radius 1 using scalars $\alpha_i>0$, producing an isomorphism of the spaces of $S$-symmetric stresses and infinitesimal motions of the original framework $(G,\op)$
 on the hyperplane and the spaces of $S^*$-symmetric stresses and infinitesimal motions of the framework $(G,q)$ on the corresponding upper half-hypersphere (hemi-hypersphere) $\SS^{d}_{+}$.
 There is no current standard form for the (orbit) rigidity matrix on the sphere that we have encountered.
  Implicitly, these frameworks are often modeled as cone frameworks with a cone joint at their center.

\begin{theorem}[Transfer between Euclidean and spherical spaces, with symmetry]
\label{thm:hemisphere}
Let $q$ be a configuration of $n$ points in $\SS^{d}_{+}$ such that the projection $\pi(q)$ from the origin (the center of the sphere) onto the hyperplane and then projected back to $\E^{d}$, is equal to $p \in \E^{dn}$. Then
\begin{enumerate}
\item[(i)]  the space of $S$-symmetric infinitesimal motions of $(G,p)$  in $\E^{d}$ is isomorphic to the space of $S^{*}$-symmetric infinitesimal motions of $(G,q)$ in $\SS_{+}^{d}$;

\item[(ii)]  the space of $S$-symmetric trivial infinitesimal motions  of $(G,p)$ in $\E^{d}$ is isomorphic to the space of $S^{*}$-symmetric trivial infinitesimal motions of $(G,q)$ in $\SS_{+}^{d}$;

\item[(iii)]   $(G,p)$ has a non-trivial $S$-symmetric  infinitesimal motion in $\E^{d}$  if and only if $(G,q)$ in $\SS_{+}^{d}$ has a non-trivial $S^{*}$-symmetric  infinitesimal motion in $\SS_{+}^{d}$;

\item[(iv)]  the space of $S$-symmetric  self-stresses  of $(G,p)$ in $\E^{d}$ is isomorphic to the space of $S^{*}$-symmetric  self-stresses of $(G,q)$ in $\SS_{+}^{d}$.

\end{enumerate}
\end{theorem}

For $S=Id$, we have given an alternative proof of the basic correspondence between spherical and Euclidean frameworks given in \cite{SalWW}.   An alternative way of describing the corresponding frameworks (without symmetry) is saying that they have the same projective coordinates.

\section{Coning finite flexes}
One of our goals in working out the details of coning and symmetry has been the opportunity coning offered for transferring finite flexes from $\E^{d}$ to $\E^{d+1}$ and to $\SS^{d}$.
   In Section \ref{subsec:transfinmot} we confirm that this transfer of finite flexes occurs when we cone with symmetry at a regular point for the class of
 $S$-symmetric realizations
(including when we simply cone at a regular point - the case of the trivial symmetry).

Section \ref{subsec:notransfinmot} includes some cautionary examples showing that at points that are not regular (for the trivial symmetry group)
 we cannot expect a transfer of finite flexes.   That is, there are cones of finitely flexing frameworks which are not finitely flexible, as well as cones which are finitely
 flexible
 even when their projections are not finitely flexible.  The larger search for more conditions that guarantee flexibility of generically infinitesimally rigid graphs
 remains an important area for continuing research.

\subsection{Transfer of finite flexes at regular points, under coning}\label{subsec:transfinmot}

We recall that one guarantee of a finite flex is that we have a non-trivial infinitesimal motion at a regular point of the graph, or a non-trivial $S$-symmetric infinitesimal motion
at an $S$-regular point of the graph (Theorem \ref{thm:flexes}).
The previous sections have confirmed the correspondence of non-trivial infinitesimal motions between cones and their projections, with and without symmetry.
  We now confirm the correspondence of regular points.

\begin{theorem}\label{thm:regptstransfer}
Let $q$ be a configuration of $n+1$ points (including the origin) in $\E^{d+1}$ such that the projection from the origin, $\pi(q)$, through the hyperplane and back to $\E^{d}$ is equal to $ p \in \E^{dn}$.
\begin{enumerate}
\item[(i)] The configuration $p\in \E^{dn}$ is a regular point of $G$ if and only if $q\in\E^{(d+1)(n+1)}$ is a regular point of $G*o$, where the cone joint is fixed at the origin;
\item[(ii)] if $(G,p)\in\mathscr{R}_{(G,S)}$ and $(G*o,q)\in\mathscr{R}_{(G*o,S^*)}$, then $p$ is an $S$-regular point of $G$ if and only if $q$ is an $S^{*}$-regular point of $G*o$, where the cone joint is fixed at the origin;
\item[(iii)] if $(G,p)\in\mathscr{R}_{(G,S)}$ and $(G*o,q)\in\mathscr{R}_{(G*o,S^*)}$, where $q\in \SS_{+}^{d}$, then $p$ is an $S$-regular point of $G$ if and only if $q$ is an $S^{*}$-regular point of $G$, where the cone joint is fixed at the origin.
\end{enumerate}
\end{theorem}
\textbf{Proof.}  (i) From the basic theorems on coning, we know that a configuration $p$ gives the maximum rank for the rigidity matrix for $G$ if and only if $q$ gives the maximum rank for the rigidity matrix for $G*o$, where the cone joint is fixed at the origin.

Moreover, moving in an open neighborhood of $p$ in $\E^{dn}$, the rank of the rigidity matrix of $G$ cannot drop immediately, but must be maintained over an open set.  The rank is determined by a maximal  minor with non-zero determinant, and small changes cannot switch this polynomial from non-zero to zero. Similarly, if we are moving within an open neighborhood of $q$ (minus the origin) in $\E^{(d+1)n}$, the rank of the cone rigidity matrix does not drop.  This confirms that $p$ is  a regular point of $G$ if and only if $q$ is a regular point of $G*o$.

(ii), (iii)  The arguments for these are entirely analogous to the previous argument.
$\square$

 While the following result is essentially a corollary to the accumulated results, we collect it as a theorem for general referencing.

\begin{theorem}[Transfer of finite flexes through coning]
\label{thm:coning flexes}
If $p$ is an $S$-regular point of $G$, and $(G*o,q)$ has symmetry $S^{*}$ with $\pi(q)=p$, then
$(G,p)$ has an $S$-symmetric finite flex if and only if $(G,q)$ has an $S^{*}$-symmetric finite flex.
\end{theorem}


\begin{remark}  [Generic implies regular]
\emph{In rigidity theory, it is common to define a configuration $p$ of $n$ points in $\mathbb{E}^d$ to be \emph{generic} if the determinant of any submatrix of the rigidity matrix $\mathbf{R}(K_n,p)$, where $K_n$ is the complete graph on $n$ vertices, is zero only if it is (identically) zero as a polynomial in the variables $p_i'$. A framework $(G,p)$ is called generic if the configuration $p$ is. Similarly, a framework $(G,p)\in \mathscr{R}_{(G,S,\Phi)}$ is defined to be \emph{$S$-generic} if any submatrix of the rigidity matrix $\mathbf{R}(K_n,p)$ is zero only if it is zero for all $p_i'$ satisfying the symmetry equations in (\ref{class}) \cite{BS1}. In other words, we obtain an $S$-generic realization of $G$ by placing the representatives of the vertex orbits under the action of $S$ into `generic' positions within their associated subspaces $U(p_i)$. Notice that it follows immediately from the definitions that the set of all generic realizations of $G$ in $\mathbb{E}^d$ !
 forms a dense open subset of all possible realizations of $G$ in $\mathbb{E}^d$, and that the set of all $S$-generic realizations of $G$ forms a dense open subset of the set $\mathscr{R}_{(G,S,\Phi)}$.
Moreover, every generic ($S$-generic) realization of $G$ is also a regular ($S$-regular) realization of $G$ in $\mathbb{E}^d$ \cite{BS6}. Thus, Theorem \ref{thm:coning flexes} applies in particular to generic and $S$-generic frameworks.
}
\end{remark}

\subsection{When finite flexes do not transfer under coning}\label{subsec:notransfinmot}
It is cautionary to recognize that the results of the previous section depend on two properties of coning:  (a) the operations of coning and projecting preserve the symmetries of the configuration, and (b) the configurations are at regular points for those symmetries.

\begin{examp}  Consider the complete bipartite graph $K_{4,4}$ in $\E^{2}$ with each of the partite sets of vertices positioned on a line (and none of the vertices lies on the point of intersection of the two lines).
This framework has a finite flex in  $\E^{2}$ if and only if the two lines are perpendicular in $\E^{2}$.  This same general condition transfers up to the sphere, as the condition that such a framework is flexible is that the two great circles are perpendicular.  Note that these properties do not require symmetry.

However, if the framework has the symmetry generated by these two lines as mirrors, the symmetry does predict the finite flexes \cite{BSWWorbit}.  So the failure of such transfer also gives insights into failures of the transfer of finite flexes induced by symmetry.
\begin{enumerate}
\item[(a)] With the mirrors (and perpendicular lines) aligned at the origin in $ \E^{2}$ the finite flexes  transfer using the prior results.

\item[(b)] If we have perpendicular lines on the sphere, and then rotate the sphere, the finite flex is preserved by the isometry. However, if we re-project with the ray through the intersection point not vertical, we can arrange that the projected image does not have perpendicular lines.  As a result, the finite flex on the sphere does not project to a finite flex in the plane.

\item[(c)] On the other hand, if we translate the framework in the plane so that the intersection point is a general point off the origin, then coning up  the framework will create a framework on the sphere without perpendicular lines.  While the infinitesimal flex will be preserved, there will be no finite flex on the sphere. \hfill $\square$
\end{enumerate} 
\end{examp}

\section{Hyperbolic space as coning in Minkowskian geometry}
\label{sec:Minkowski}
We can extend the processes of the previous sections to a process which will transfer results for Euclidean frameworks in $\E^{d}$ through to the corresponding frameworks in hyperbolic space $\H^{d}$. This process involves embedding $\E^{d}$ into a hyperplane in Minkowskian geometry, and then coning to the origin so that this cone framework lies on a hyper-hyperboloid which carries the metric  of hyperbolic geometry.

Again the entire process extends to  symmetry groups, provided we are careful about the orientation of the plane we use in the Minkowskian geometry. We begin with the basic structure of rigidity in Minkowskian space.

\subsection{Minkowskian frameworks}
In this section we will focus on frameworks embedded in the $(d+1)$-dimensional Minkowskian space $\mathbb{M}^{d+1}$. This will follow precisely the steps in \S2.1.  However, since the rigidity theory in Minkowskian space is not widely studied \cite{alex}, we state the basic definitions and results.

As a basic geometry, the $(d+1)$-dimensional Minkowskian space $\mathbb{M}^{d+1}$ is the $(d+1)$-dimensional metric space with the metric
     $$\|(a_{1}, \ldots, a_{d}, a_{d+1})\|_{\mathbb{M}}^{2} = a_{1}^{2}+\dots+a_{d}^{2} - a_{d+1}^{2}.$$
A bar and joint framework in Minkowskian $(d+1)$-space is a graph $G$, with vertices $V(G)$ and edges $E(G)$,  and a map $p:V(G)\rightarrow \M^{d+1}$ such that for $i,j\in V(G)$, we have $p_{i}\neq p_{j}$ (the joints are distinct).

This slight change in the metric means that we make a slight change in the equations defining an infinitesimal motion and the corresponding rigidity matrix.   For a fixed ordering of the edges of a graph $G$ with vertex set $V(G)=\{1,\ldots,n\}$, we define the \emph{edge function} $f^{\M}_{G}:\mathbb{M}^{(d+1)n}\to \mathbb{R}^{|E(G)|}$  by
\begin{displaymath}
f^{\M}_{G}\big(p_{1},\ldots,p_{n}\big)=\big(\ldots, \|p_{i}-p_{j}\|_{\M}^2,\ldots\big)\textrm{, }
\end{displaymath}
where $\{i,j\}\in E(G)$, $p_{i}:=p(i)\in \mathbb{M}^{d+1}$ for all $i=1,\ldots,n$, and $\| \cdot \|_{\M}$ denotes the Minkowskian norm in $\mathbb{\M}^{d+1}$.

If $(G,p)$ is a $(d+1)$-dimensional framework with $n$ vertices, then $(f^{\M}_{G})^{-1}\big(f^{\M}_{G}(p)\big)$ is the set of all configurations $q$ of $n$ points in $\mathbb{M}^{d+1}$ with the property that corresponding bars of the frameworks $(G,p)$ and $(G,q)$ have the same length. In particular, we clearly have
$(f^{\M})^{-1}_{K_{n}}\big(f^{\M}_{K_{n}}(p)\big)\subseteq (f^{\M}_{G})^{-1}\big(f^{\M}_{G}(p)\big)$, where $K_{n}$ is the complete graph on $V(G)$.

An analytic path $x:[0,1]\to \mathbb{M}^{(d+1)n}$ is called a \emph{finite motion} of $(G,p)$ if $x(0)=p$ and $x(t)\in (f^{\M}_{G})^{-1}\big(f^{\M}_{G}(p)\big)$ for all
$t\in [0,1]$. Further, $x$ is called a \emph{finite rigid motion} (or \emph{trivial finite motion}) if $x(t)\in (f^{\M}_{K_{n}})^{-1}\big(f^{\M}_{K_{n}}(p)\big)$ for all $t\in [0,1]$, and $x$ is called a \emph{finite flex} (or \emph{non-trivial finite motion}) of $(G,p)$ if $x(t)\notin (f^{\M}_{K_{n}})^{-1}\big(f^{\M}_{K_{n}}(p)\big)$ for all $t\in (0,1]$.  While all of the translations are trivial motions (as in $\E^{d+1}$) the `rotations' take a different form, moving a point along a hyperbolic path of `constant distance' rather than a circular path, because of the modified metric.

We say that $(G,p)$ is \emph{rigid in $\M^{d+1}$} if every finite motion of $(G,p)$ is trivial; otherwise $(G,p)$ is called \emph{flexible in $\M^{d+1}$}.

Given a framework $(G,p)$ in $\M^{d+1}$, the
\emph{rigidity matrix} of $(G,p)$ is the $|E(G)| \times (d+1)n$ matrix $\mathbf{R^{\M}}(G,p)=\frac{1}{2}df^{\M}_G(p)$, where $df^{\M}_G(p)$ denotes the Jacobian matrix of the edge function $f^{\M}_G$, evaluated at the point $p$.  However, if we are more explicit, we see that the entries are not quite $(p_{i} - p_{j})$.  Rather if we write $\hat{a} = (a_{1},\ldots, a_{d}, - a_{d+1})$, they have the form:
$\widehat{(p_{i} - p_{j})} = (\hat p_{i} - \hat p_{j})$.
So the Minkowskian rigidity matrix for the framework $(G,p)$ in $\M^{d+1}$ has the form
\begin{displaymath} \mathbf{R}^{\M}(G,p)=\bordermatrix{& & & & i & & & & j & & & \cr & & & &  & & \vdots & &  & & &
\cr \{i,j\} & 0 & \ldots &  0 & (\hat p_{i}-\hat p_{j}) & 0 & \ldots & 0 & (\hat p_{j}- \hat p_{i}) &  0 &  \ldots&  0 \cr & & & &  & & \vdots & &  & & &
}
\textrm{.}\end{displaymath}
The solutions to the equation $ \mathbf{R}^{\M}(G,p) u = 0$ are the infinitesimal motions in Minkowskian space.   Provided the joints of the framework affinely span the space $\M^{d+1}$, the trivial infinitesimal motions are the solutions to the equation: $ \mathbf{R}^{\M}(K_{n},p) u$. It remains true that this space is of dimension $\binom{(d+1)+1}{2}=\binom{d+2}{2}$.

We say that $(G,p)$ is \emph{infinitesimally rigid in $\M^{d+1}$} if every infinitesimal motion of $(G,p)$ is an infinitesimal rigid motion. Otherwise $(G,p)$ is said to be \emph{infinitesimally flexible}.

Clearly, if the joints of $(G,p)$ affinely span all of $\mathbb{M}^{d+1}$, then $\textrm{nullity }\big(\mathbf{R}^{\M}(G,p)\big)\geq \binom{d+2}{2}$, and $(G,p)$ is infinitesimally rigid in ${\M}^{d+1}$ if and only if $\textrm{nullity } \big(\mathbf{R}^{\M}(G,p)\big)=\binom{d+2}{2}$ or equivalently, $\textrm{rank }\big(\mathbf{R}^{\M}(G,p)\big)=(d+1)n - \binom{d+2}{2}$.

An infinitesimally rigid framework is always rigid. As for Euclidean space, the converse, however, does not hold in general.

A \emph{self-stress} of a framework $(G,p)$ with $V(G)=\{1, \ldots, n\}$ remains a row dependence of the modified rigidity matrix.  That makes it a  function  $\omega:E(G)\to \mathbb{R}$ such that at each joint $p_i$ of $(G,p)$ we have
\begin{displaymath}
\sum_{j :\{i,j\}\in E(G)}\omega_{ij}(\hat p_{i}- \hat p_{j})=0 \textrm{,}
\end{displaymath}
where $\omega_{ij}$ denotes $\omega(\{i,j\})$ for all $\{i,j\}\in E(G)$. Note that if we identify a self-stress $\omega$ with a row vector in $\mathbb{R}^{|E(G)|}$ (by using the order on $E(G)$), then we have $\omega\mathbf{R}^{\M}(G,p)=0$.

If $(G,p)$ has a non-zero self-stress, then $(G,p)$ is said to be \emph{dependent in $\M^{d+1}$} (since in this case there exists a linear dependency among the row vectors of $\mathbf{R}^{\M}(G,p)$). Otherwise, $(G,p)$ is  \emph{independent in $\M^{d+1}$ }.  A framework which is both independent and infinitesimally rigid is called \emph{isostatic in $\M^{d+1}$ }.

A configuration $p$ of $n$ points in $\mathbb{M}^{d+1}$ is called a \emph{regular point of the graph $G$}, if $\textrm{\rm rank\,}\big(\mathbf{R}^{\M}(G,p)\big)\geq\textrm{\rm rank\,} \big(\mathbf{R}^{\M}(G,q)\big)$ for all $q\in \mathbb{M}^{(d+1)n}$.
A framework $(G,p)$ is said to be \emph{regular} if $p$ is a regular point of $G$.

It follows immediately from this definition that the set of all regular realizations of a graph $G$ in $\mathbb{M}^{d+1}$ forms a dense open subset of all possible realizations of $G$ in $\mathbb{M}^{d+1}$. Moreover, note that the infinitesimal rigidity of a regular realization of $G$ depends only on the underlying graph $G$ and not on the particular realization.

The results of Asimov and Roth in \cite{asiroth} generalize in a natural way so that for regular frameworks in  $\M^{d+1}$, infinitesimal rigidity and rigidity are  equivalent. This result continues to provide a key tool for detecting finite flexes in frameworks.

\subsection{Coning from Euclidean space into Minkowskian space}

In this section, we will embed the Euclidean space $\E^{d}$ as the hyperplane $x_{d+1} = 1$ in the Minkowskian space $\M^{d+1}$.  Notice that the metric among points within this hyperplane is the same as the metric in the Euclidean space, since the final coordinate cancels out. Given a framework $(G,p)$ in $\E^d$, we embed $(G,p)$ into the hyperplane $x_{d+1}=1$ of $\M^{d+1}$ via $\op_{i}=(p_{i},1)\in \M^{d+1}$, and then cone the resulting framework to the origin $o$ in $\M^{d+1}$ with $n$ new edges, adding $n$ rows to the new rigidity matrix in $\M^{d+1}$, creating  the {\emph{cone graph, $G*o$,} and the \emph{cone framework} $(G*o,\ps)$}. This procedure follows the basic steps of Section \ref{sec:conewithoutsym}.


The
\emph{cone rigidity matrix in $\M^{d+1}$} of the cone framework $(G*o,\ps)$ is the $(|E(G)|+n) \times (d+1)n$ matrix
\begin{displaymath}{\RC}^{\M}(G*o,\ps)=\bordermatrix{& & & & i & & & & j & & & \cr & & & &  & & \vdots & &  & & &
\cr
 \{i,j\} & 0 & \ldots &  0 & (\hat{\op}_{i}-\hat{\op}_{j}) & 0 & \ldots & 0 & (\hat{\op}_{j}-\hat{\op}_{i}) &  0 &  \ldots&  0
 \cr & & & &  & & \vdots & &  & & &
\cr \{0,i\} & 0 & \ldots &  0 &\hat{\op}_{i}& 0 & \ldots & 0 & 0&  0 &  \ldots&  0
\cr & & & &  & & \vdots & &  & & &
\cr  \{0,j\} & 0 & \ldots &  0 & 0& 0 & \ldots & 0 & \hat{\op}_{j} &  0 &  \ldots&  0
\cr & & & &  & & \vdots & &  & & &
}
\textrm{.}\end{displaymath}
As before, we have added $n$ rows and $n$ columns, and $(\hat {\op}_{i}-\hat{\op}_{j})$ is zero in each added column. Moreover, for each added column (under vertex $i$) there is exactly one added row which is non-zero in this column: $\hat{\op}_{i}$ has a $-1$ in this column.  Thus we have increased the rank by $n$, and preserved the dimension of the kernel.  For convenience, we have given a modified cone rigidity matrix which omits columns for the cone joint (the origin), and the infinitesimal motions are restricted to those which fix the origin.   This takes infinitesimal motions of $(G,p)$ in $\E^{d}$ to infinitesimal motions of $(G*o,\ps)$ in $\M^{d+1}$.

When we apply this same process to a complete graph on the vertices, we see that the trivial infinitesimal motions in $\E^{d}$ go to trivial infinitesimal motions in $\M^{d+1}$ which fix the origin.

\begin{figure}
    \begin{center}
  \subfloat[] {\label{fig:conevelocity}{\includegraphics [width=.27\textwidth]{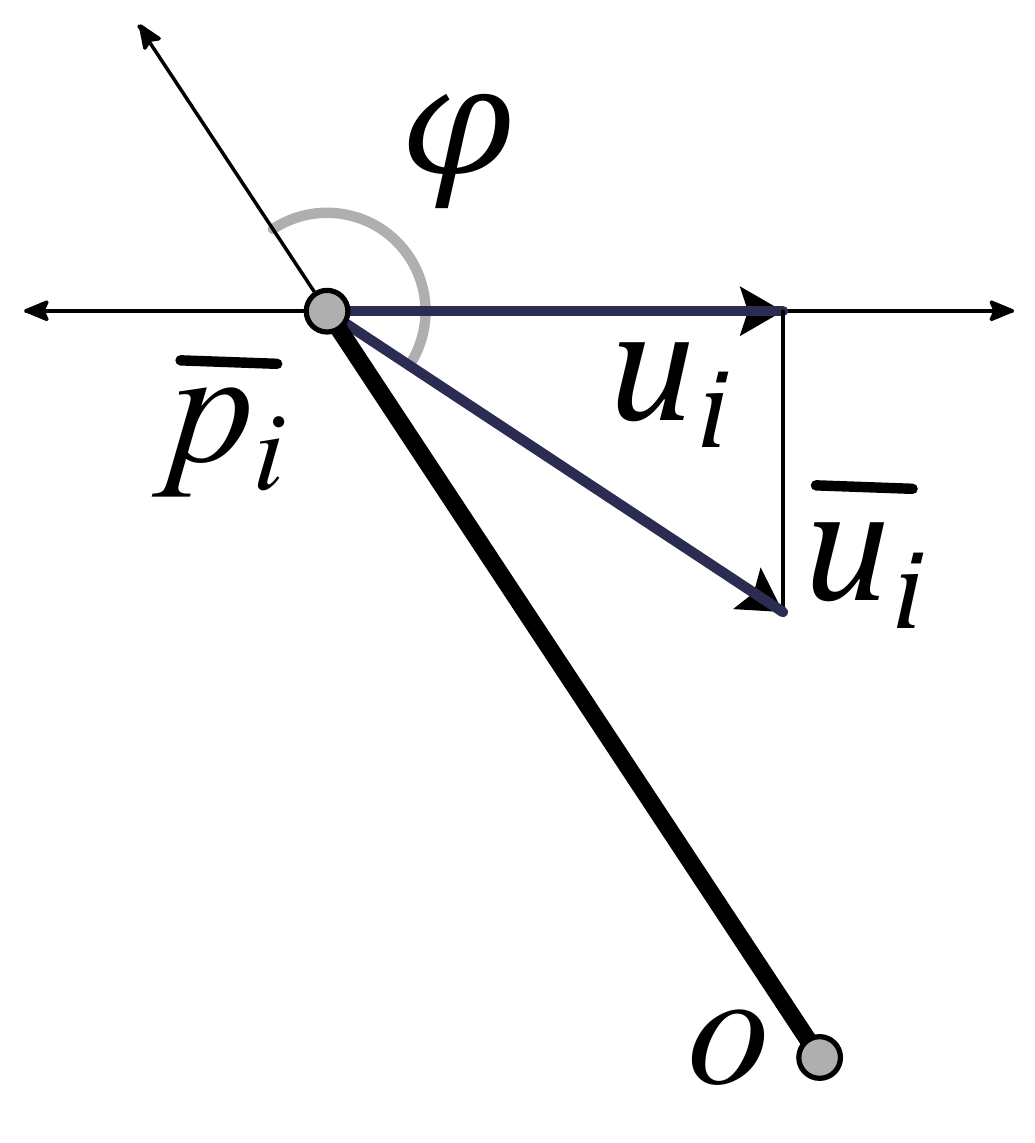}}}\quad
 \subfloat[]{\label{fig:hyperbolavelocity1}{\includegraphics[width=.33 \textwidth]{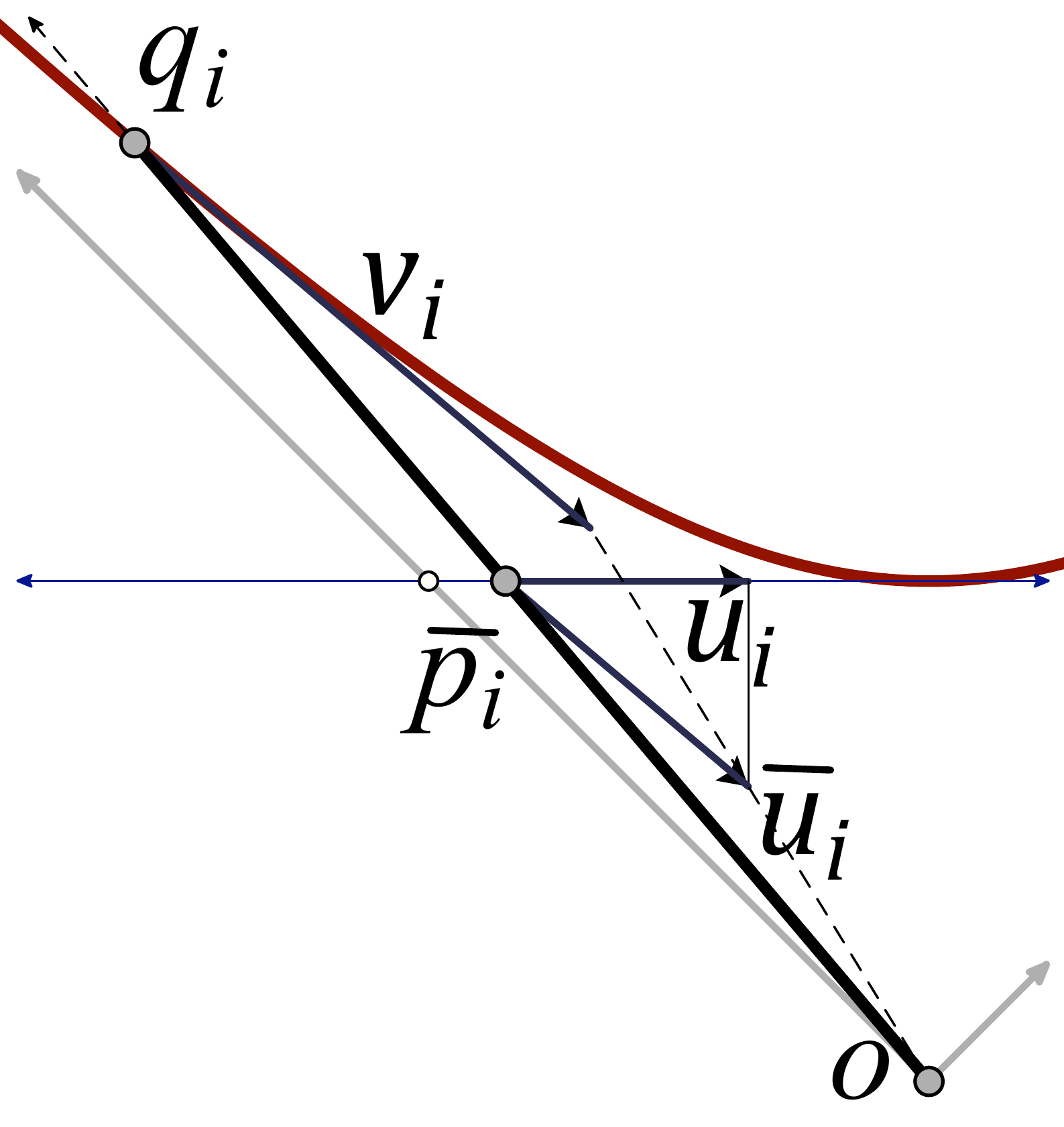}}}  \quad
  \subfloat[]{\label{fig:hyperbolavelocity2}{\includegraphics[width=.33 \textwidth]{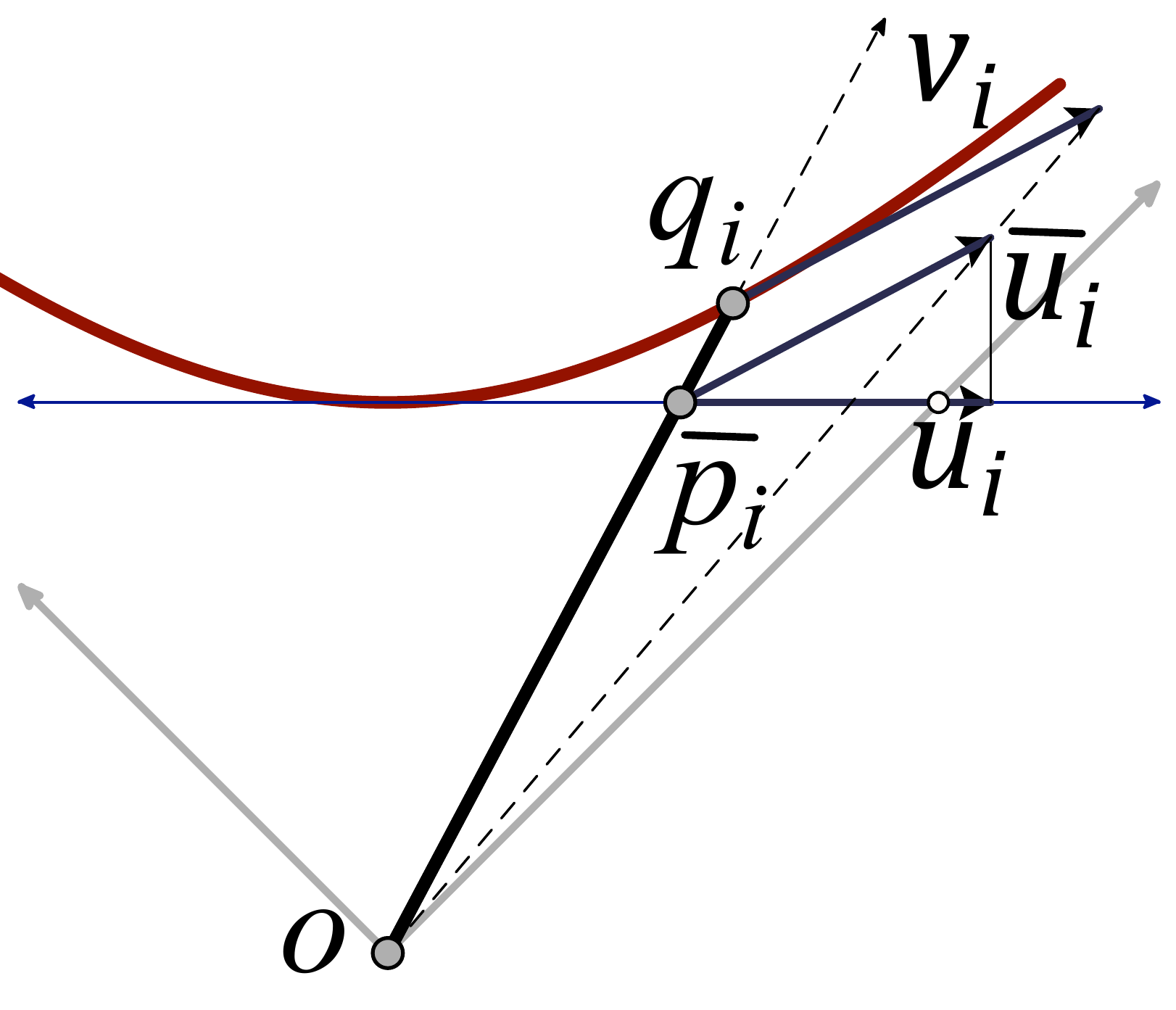}} }
      \end{center}
  \caption{Moving a velocity under coning into the Minkowskian space requires that we use the Minkowskian measure of `right angle' ($\phi=90^{}$) (a). We then push the vertices (and velocities) onto the hyperboloid (c), (d). If $u_i=u_j$ then the new velocities have the same length in the Minkowskian metric ($||v_i||_\M =||v_j||_\M$). }
    \label{fig:MinkowskiCone}
    \end{figure}

What does the transfer of velocities look like, in practice?  Consider Figure~\ref{fig:MinkowskiCone}(a) which shows the process from the line to the plane.
We keep the same first $d$ coordinates of the velocities, and add whatever last entry will make the vector perpendicular to the bar from the origin to the joint, in the measure of angle in Minkowskian space.  That is, we take the unique vector which is perpendicular to the coning bar and projects orthogonally onto the previous velocity.

\subsection{Coning with symmetry in $\M^{d+1}$}

This subsection will follow the steps in Section \ref{sec:conewithsym}.
As before, we first show that for a given framework $(G,p)$ in $\mathscr{R}_{(G,S)}$ in $\E^{d}$, the cone framework $(G*o,\ps)$ in $\M^{d+1}$ has the equivalent  symmetry group $S^{*}$, lying in the corresponding space $\mathscr{R}^\M_{(G,S^*)}$.  It also has the same symmetric infinitesimal rigidity properties as $(G,p)$.

\begin{defin} \emph{Let $S$ be a symmetry group in dimension $d$ and let $x\in S$. Recall that we denote the matrix which represents $x$ w.r.t. the canonical basis of $\mathbb{E}^d$ by $M_x$. We define $x^*$ to be the transformation of $\mathbb{M}^{d+1}$ which is represented by the matrix $$M_{x^*}=\left( \begin{array}{cccc} & & & 0\\& M_x & & \vdots\\ & & & 0\\0 & \cdots & 0 &1 \end{array} \right)$$
w.r.t. the canonical basis of $\mathbb{M}^{d+1}$.
 This is still an isometry in $\M^{d+1}$, since
\begin{eqnarray*} \|((a_{1},\ldots,a_{d},a_{d+1})M_{x^{*}}^T)\|_{\M}^2 & = & \|(a M_{x}^T,a_{d+1})\|_{\M}^2\\ & = & \|\pi(a)M_{x}^T\|_{\E}^2- a_{d+1}^2\\ & = & \|\pi(a)\|_{\E}^2- a_{d+1}^2 \\ & = & \|(a_{1},\ldots,a_{d},a_{d+1})\|_{\M}^2,\end{eqnarray*}
where $a=(a_1,\ldots,a_d)$ and $\pi(a_1,\ldots,a_d,a_{d+1}):=(a_1,\ldots,a_d)$.
We let $S^*$ be the group in dimension $(d+1)$ which has the elements $\{x^*:x\in S\}$.}\end{defin}

We immediately have the corresponding version of Theorem~\ref{thm:groupextend}
\begin{theorem}
A framework $(G,p)$ is an element of $\mathscr{R}_{(G,S)}$ if and only if the cone framework $(G*o,\ps)$ is an element of $\mathscr{R}^\M_{(G*o,S^*)}$.
\end{theorem}

We now can present the complete analog of the previous results, for coning from $\E^{d}$ into $\M^{d+1}$, provided that  $\E^{d}$  is placed as a hyperplane orthogonal to the one negative component in the signature of the metric.

\begin{theorem}
Given a $d$-dimensional Euclidean framework $(G,p)$ in $\mathscr{R}_{(G,S)}$ and the corresponding $(d+1)$-dimensional Minkowskian cone framework $(G*o,\ps)$ in $\mathscr{R}^\M_{(G,S^*)}$, we have that
\begin{enumerate}
\item[(i)]  the space of $S$-symmetric infinitesimal motions of $(G,p)$ in $\E^{d}$  is isomorphic to the space of  $S^*$-symmetric infinitesimal motions of $(G*o,\ps)$ in $\M^{d+1}$, with the cone joint of $(G*o,\ps)$ fixed at the origin;
\item[(ii)]  the space of trivial $S$-symmetric infinitesimal motions  of $(G,p)$ in $\E^{d}$ is isomorphic to the space of  $S^*$-symmetric trivial infinitesimal motions of $(G*o,\ps)$ in $\M^{d+1}$, with the cone joint of $(G*o,\ps)$ fixed at the origin;
\item[(iii)]   $(G,p)$ has a non-trivial $S$-symmetric infinitesimal motion in $\E^{d}$  if and only if $(G*o,\ps)$  has a non-trivial $S^*$-symmetric infinitesimal motion in $\M^{d+1}$;
\item[(iv)]  the space of $S$-symmetric self-stresses of $(G,p)$  in $\E^{d}$ is isomorphic to the space of $S^*$-symmetric self-stresses of $(G*o,\ps)$ in $\M^{d+1}$.
\end{enumerate}
\end{theorem}
\textbf{Proof.} We use the same procedure as in the proofs of Theorems~\ref{thm:ConingGeneral} and \ref{th:rigpropcone}.
As before, for a framework $(G,p)\in\mathscr{R}_{(G,S)}$ in $\E^{d}$ whose vertex orbits under the action of $S$ are represented by the vertices $1\ldots, k$, the orbit rigidity matrix  $\mathbf{O}(G,p,S)$ of $(G,p)$ is the $r \times c$ matrix
\begin{displaymath}\renewcommand{\arraystretch}{0.7}
    \bordermatrix{  & &  i &  &  j & \cr &&& \vdots &&\cr \{i,x(j)\} & 0 \ldots 0 & (p_i-x(p_j))\mathbf{M}_{i} & 0  \ldots  0 & (p_j-x^{-1}(p_i))\mathbf{M}_{j} & 0  \ldots  0\cr
    &&& \vdots &&\cr
    \{i,y(i)\} & 0  \ldots  0 & (2p_i-y(p_i)-y^{-1}(p_i))\mathbf{M}_{i} & 0  \ldots  0 &0 & 0  \ldots  0\cr
      &&& \vdots &&
    }\textrm{,}
    \end{displaymath}
    where $c=\sum_{i=1}^{k} c_i$.

We embed the framework $(G,p)$ into the hyperplane $x_{d+1}=1$ of $\mathbb{M}^{d+1}$ via $\op_{i}=(p_{i},1) \in \M^{d+1}$.  We then cone this to the origin $o$ in $\M^{d+1}$. This gives rise to $k$ new rows in the orbit rigidity matrix $\mathbf{\OC}^{\M}(G*o,\ps,S^*)$ of the framework $(G*o,\ps)\in\mathscr{R}_{(G*o,S^*)}$.

For each $i=1,\ldots, k$, we let $\mathscr{B}^*_{i}$ be the basis of the space $U(\ps_i)$ which consists of the basis vectors $\{(w,0): w\in \mathscr{B}_{i}\}$ and the additional basis vector $e_{d+1}$ - the $(d+1)$st canonical basis vector of $\mathbb{M}^{d+1}$. Further, we let $\mathbf{M}_i^*$  be the matrix whose columns are the coordinate vectors of  $\mathscr{B}^*_{i}$ relative to the canonical basis of $\mathbb{M}^{d+1}$; that is, $$\mathbf{M}_i^*=\left( \begin{array}{cccc} & & & 0\\& \mathbf{M}_i & & \vdots\\ & & & 0\\0 & \cdots & 0 &1 \end{array} \right).$$
The orbit rigidity matrix $\mathbf{\OC}^{\M}(G*o,\ps,S^*)$  has the form
   $$\renewcommand{\arraystretch}{0.7}
    \bordermatrix{  & &  i &  &  j & \cr &&& \vdots &&\cr \{i,x(j)\} & 0 \ldots 0 & (\hat \op_i-x^*(\hat \op_j))\mathbf{M}_i^* & 0  \ldots  0 & (\hat \op_j-x^{*-1}(\hat \op_i))\mathbf{M}_j^* & 0  \ldots  0\cr
    &&& \vdots &&\cr
     \{i,y(i)\}& 0  \ldots  0 & (2 \hat \op_i-y^*(\hat \op_i)-y^{*-1}(\hat \op_i))\mathbf{M}_i^* & 0  \ldots  0 &0 & 0  \ldots  0\cr
      &&& \vdots &&
      \cr \{0,i\}& 0 \ldots  0 &\hat{\op}_{i}\mathbf{M}_i^*& 0 \ldots 0 & 0&  0  \ldots  0
\cr & & &  \vdots  & &
\cr \{0,j\}  & 0 \ldots  0 & 0& 0  \ldots  0 & \hat{\op}_{j}\mathbf{M}_j^* &  0  \ldots  0
\cr & & &  \vdots  & &
      }\textrm{.}
    $$

Note that we have added $k$ rows and $k$ columns, and that for every edge $\{i,x(j)\}$, $i\neq j$, and every edge $\{i,y(i)\}$ of $G$, the vectors $(\hat \op_i-x^*(\hat \op_j))\mathbf{M}_i$ and $(2 \hat \op_i-y^*(\hat \op_i)-y^{*-1}(\hat \op_i))\mathbf{M}_i^*$ are zero in the added column.  Moreover, for each added column (under vertex $i$) there is exactly one added row which is non-zero in this column: $\hat{\op}_{i}\mathbf{M}_i^*$ has a $-1$ in this column.  Thus we have increased the rank by $k$, and preserved the dimension of the kernel.

If we apply this process to the complete graph on the vertices, we see that the dimension $m$ of the space of trivial $S$-symmetric infinitesimal motions of $(G,p)$ is equal to the dimension $m^*$ of the space of trivial $S^*$-symmetric infinitesimal motions of $(G*o,\ps)$ (with the cone joint fixed at the origin). This gives the result. $\square$

\subsection{Pushing vertex orbits in Minkowskian space}

Following the same steps as used in \S3.2-\S3.5, we can push and pull orbits of joints along the cone rays to create the framework $(G*o,q,S^{*})$.   All of the steps involve simple row reductions in the orbit matrix $\mathbf{\OC}^{\M}(G*o,\ps,S^*)$ plus multiplying the rows and columns by non-zero constants.  Without repeating all the details, we can conclude with the analog of Theorem~\ref{thmtranconwithsym}.

\begin{theorem} [Symmetry coning from Euclidean into Minkowskian space] Given a configuration $q$ of $n+1$ points (including the origin) in $\M^{d+1}$ such that the projection $\pi(q)$ from the origin onto the hyperplane and then back to $\E^d$, is equal to $p \in \E^{dn}$, we have:
\begin{enumerate}
\item[(i)]  the space of $S$-symmetric infinitesimal motions of $(G,p)$ in $\E^{d}$ is isomorphic to the space of $S^{*}$-symmetric infinitesimal motions of $(G*o,q)$ in $\M^{d+1}$, with the cone joint of $(G*o,q)$ fixed at the origin;
\item[(ii)]  the space of $S$-symmetric trivial infinitesimal motions of $(G,p)$ in $\E^{d}$ is isomorphic to the space of $S^{*}$-symmetric trivial infinitesimal motions of $(G*o,q)$ in $\M^{d+1}$, with the cone joint of $(G*o,q)$ fixed at the origin;
\item[(iii)]   $(G,p)$ has a non-trivial $S$-symmetric  infinitesimal motion in $\E^{d}$  if and only if $(G*o,q)$ has a non-trivial $S^{*}$-symmetric  infinitesimal motion in $\M^{d+1}$;
\item[(iv)]  the space of $S$-symmetric  self-stresses of $(G,p)$ in $\E^{d}$ is isomorphic to the space of $S^{*}$-symmetric  self-stresses of $(G*o,q)$ in $\M^{d+1}$.
\end{enumerate}
\end{theorem}

\subsection{Pushing vertex orbits onto the hyper-hyperboloid}

\begin{figure}
    \begin{center}
 \subfloat[]{\includegraphics[width=2.4in]{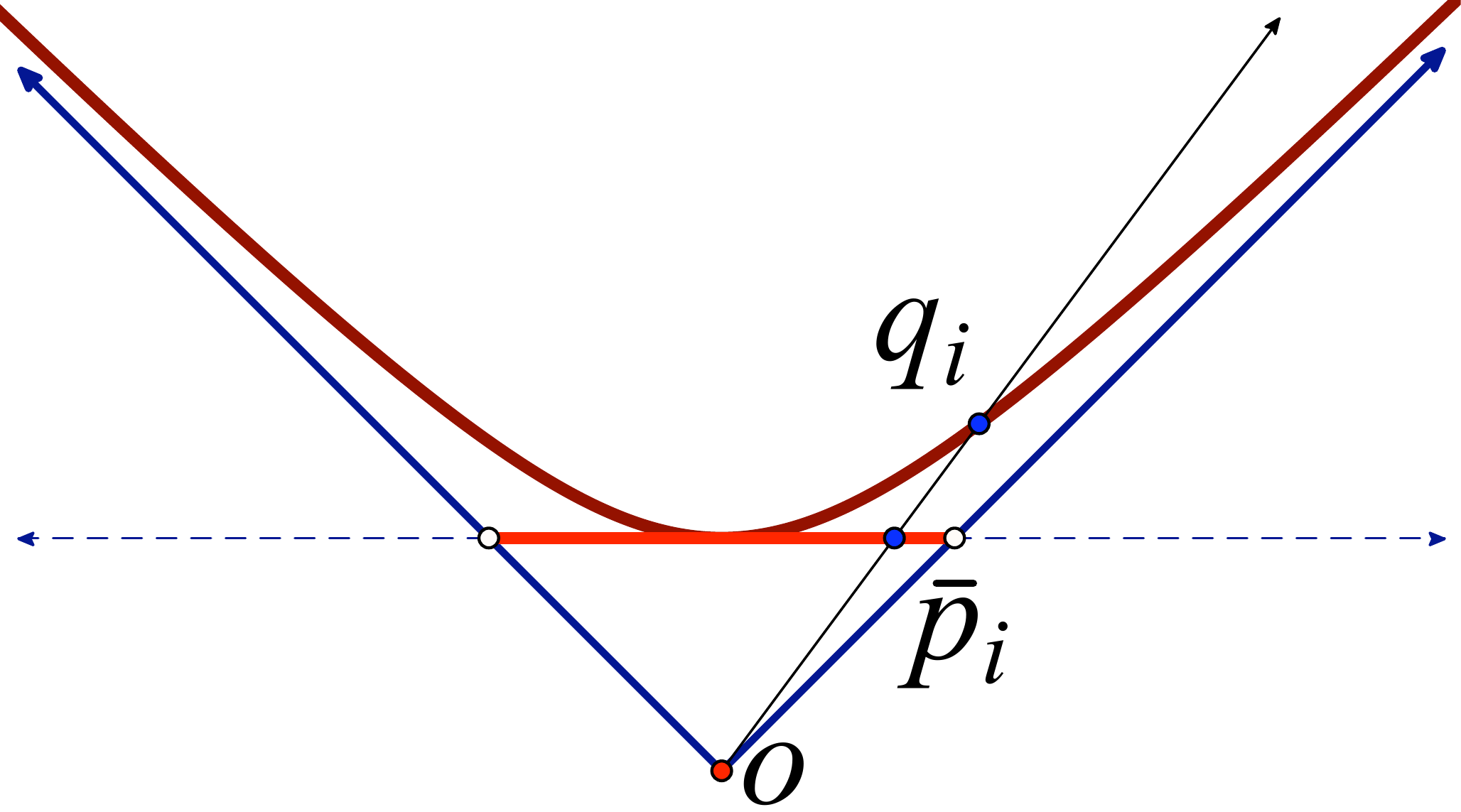}}  \quad
  \subfloat[]{\label{fig:conevelocity}\includegraphics[width=1.7in]{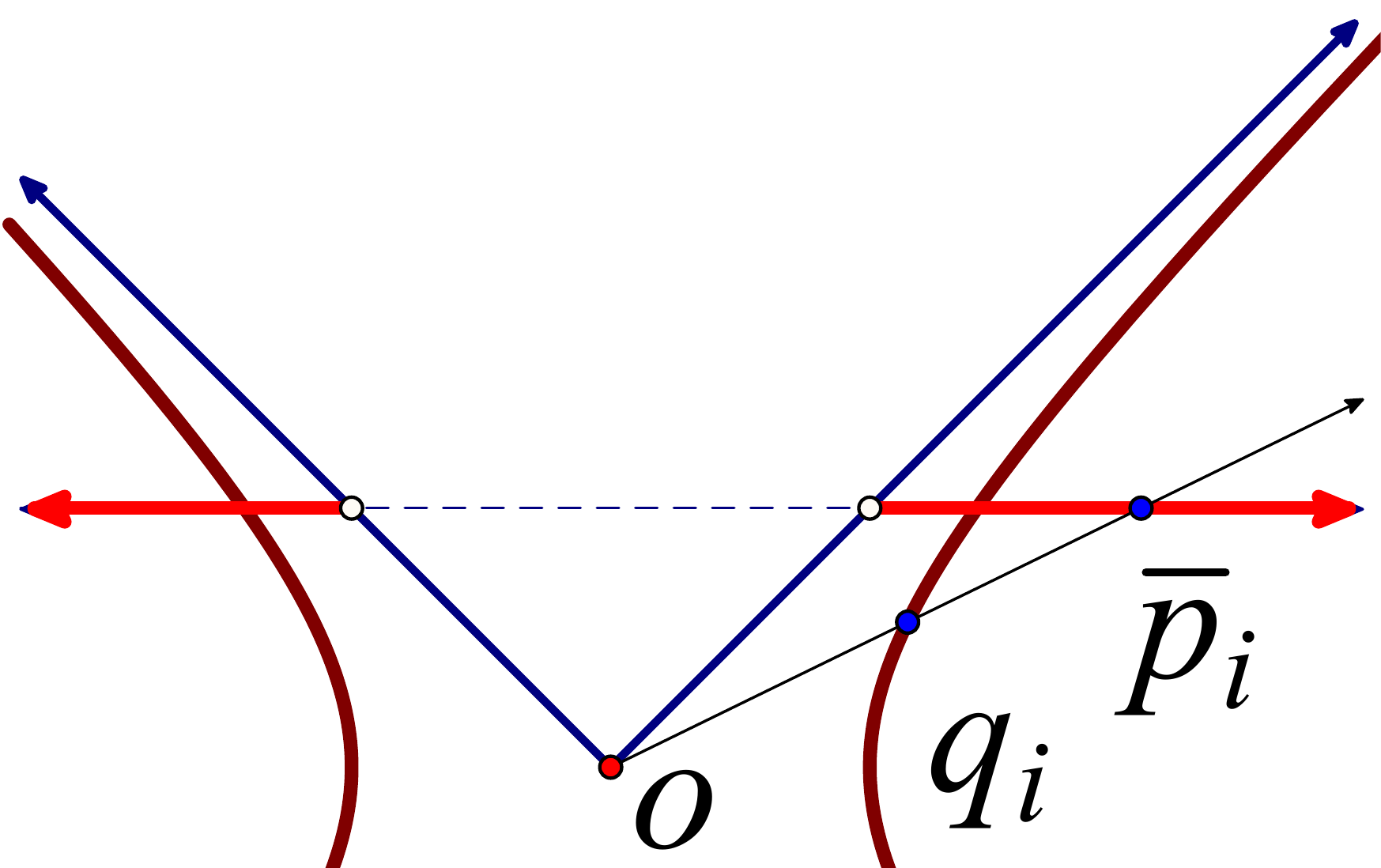}}   \quad\quad\quad
      \end{center}
  \caption{The Minkowskian sphere of radius squared equal to -1 projects to an interval on the line  (a) and the Minkowskian sphere of radius squared equal to 1 projects to the exterior of that interval (b).}
    \label{fig:MinkowskiSpheres}
    \end{figure}

By carefully selecting the scalars $\alpha_{i}$, we can ensure that all orbits of joints are placed on the hyper-hyperboloid
$\|a\|^{2}_{\M} = a_{1}^{2} + \ldots + a_{d}^{2} - a_{d+1}^{2} = -1$, provided that the original points in $\E^{d}$ satisfy
$\|\pi(a)\|_{\E}  < 1$ (Figure~\ref{fig:MinkowskiSpheres}(a)).     When we achieve this placement, the metric within the hyper-hyperpoloid is the hyperbolic
metric of the space $\H^{d}$.   Moreover, the symmetries of $S$ remain as $S^{*}$ - symmetries of the hyper-hyperboloid, or symmetries within the hyperbolic space $\H^{d}$.  Together, this translation and the basic coning results give the following equivalence.

\begin{theorem}
\label{thm:hyperbolic} Given a configuration $q$ of $n$ points in $\H^{d}$ such that the projection $\pi(q)$ from the origin onto the hyperplane and then back to $\E^d$ is equal to $p\in \E^{dn}$, we have:
\begin{enumerate}
\item[(i)]   the space of $S$-symmetric infinitesimal motions of $(G,p)$ in $\E^{d}$ is isomorphic to the space of $S^{*}$-symmetric infinitesimal motions of $(G,q)$ in $\H^{d}$;
\item[(ii)]   the space of $S$-symmetric trivial infinitesimal motions  of $(G,p)$ in $\E^{d}$ is isomorphic to the space of $S^{*}$-symmetric trivial infinitesimal motions of $(G,q)$ in $\H^{d}$;
\item[(iii)]    $(G,p)$ has a non-trivial $S$-symmetric  infinitesimal motion in $\E^{d}$  if and only if $(G,q)$ has a non-trivial $S^{*}$-symmetric  infinitesimal motion in $\H^{d}$;
\item[(iv)]   the space of $S$-symmetric  self-stresses  of $(G,p)$ in $\E^{d}$ is isomorphic to the space of $S^{*}$-symmetric  self-stresses of $(G,q)$ in $\H^{d}$.
\end{enumerate}
\end{theorem}

Note that for any framework in $\E^{d}$ it is just a matter of dilating the framework towards the origin to ensure that the framework lives within the disc $\|a\|_{\E} <1$.

There is a second geometry which lives on the alternative hyper-hyperboloid
$\|a\|^{2}_{\M} = a_{1}^{2} + \ldots + a_{d}^{2} - a_{d+1}^{2} = 1$ (Figure~\ref{fig:MinkowskiSpheres}(b)).  This is called the de Sitter geometry ${\DS}^{d}$ \cite{Coxeter}.
If the original framework $(G,p)$ has no fixed joints, then we can instead expand $(G,p)$ so that the entire set of joints and bars lives outside of the disc $\|a\|_{\E} \leq 1$. We can then push the expanded framework onto the de Sitter space.
All of the results for pushing from cones shift to this second surface without any extra effort, including Theorem~\ref{thm:hyperbolic}.

\subsection{Finite flexes in hyperbolic space}

It is now clear that the methods and results of \S4.1 for the transfer of finite flexes at regular points between $\E^{d}$ and $\SS^{d}$ immediately extend to a transfer of finite flexes at regular points between $\E^{d}$  and $\H^{d}$.  We do not repeat the arguments, but simply state the final conclusion.

\begin{theorem}[Transfer of finite flexes through coning]
\label{thm:coninghyperbolicflexes}
If $p$ is an $S$-regular point of $G$ in $\E^{dn}$, and $(G,q)$ in $\H^{d}$ has symmetry $S^{*}$ with $\pi(q)=p$, then
$(G,p)$ has an $S$-symmetric finite flex in $\E^{d}$ if and only if $(G,q)$ has an $S^{*}$-symmetric finite flex in $\H^{d}$.
\end{theorem}

\section{Further discussion}
\subsection {Transfer of body-bar flexibility between spaces}

Recently, there has been an increasing attention to a special class of frameworks, the body-bar frameworks, including work on symmetries of body-bar frameworks \cite{WWbb, gsw}.   The results presented here can be used to provide a transfer of infinitesimal rigidity results  for these frameworks between $\E^{d}$, $\SS_{+}^{d}$, and $\H^{d}$.   The idea is to record a body-bar framework as a special form of a bar-and-joint framework.

Let $H$ be a multigraph with minimum degree at least one.
The {\it body-bar graph induced by
$H$}, denoted by $G_H$, is the graph obtained from $H$ by
replacing each vertex $v\in V(H)$ by a complete graph $B_v$ (a `body') on
$deg_H(v)$ vertices and replacing each edge $\{u,v\}$ of $H$ by an edge (a `bar')
between $B_u$ and $B_v$ in such a way that the bars are pairwise
disjoint.

Notice that this is a combinatorial correspondence, so the graph $G_{H}$ can be built as a framework in $\E^{d}$, $\SS_{+}^{d}$, or
$\H^{d}$.   Since these are bar-and-joint frameworks, they can have symmetries, and can be coned in exactly the way we have discussed throughout this paper.  This permits the full transfer of infinitesimal rigidity results among the metrics.   The coned framework is not a body-bar framework, but the cone vertex disappears in the explicit graph when we move on to the frameworks in
$\SS_{+}^{d}$ or $\H^{d}$.
Moreover, if we have specific configurations for the vertices of the bars, we can define regular realizations, and transfer finite motions, with and without symmetry.

The one caution is that in the setting of infinitesimal and generic rigidity, we often consider two body-bar frameworks $(G_{H}, p)$ and $(G_{H}, p')$  equivalent provided that the lines of the edges are the same in the two frameworks.  Such equivalence does not necessarily preserve all of the symmetry, nor does it necessarily preserve finite motions.   For this reason we state the transferred results for explicit choices of the configuration of the joints on the ends of the bars.

\begin{theorem} Given a body-bar framework $(G_{H},p)$ with symmetry group $S$ in $\E^{d}$, a body-bar framework $(G_{H},q)$ with the corresponding symmetry group $S^{*}$ in $\SS_{+}^{d}$, where $\pi(q)=p$, and a body-bar  framework $(G_{H},r)$ with symmetry group $S^{*}$ in $\H^{d}$, where $\pi(r)=p$, the following are equivalent:
\begin{enumerate}
\item[(i)]  $(G_{H},p)$ has a non-trivial $S$-symmetric infinitesimal motion in $\E^{d}$;
\item[(ii)]  $(G_{H},q)$  has a non-trivial $S^*$-symmetric infinitesimal motion in $\SS_{+}^{d}$;
\item[(iii)]  $(G_{H},r)$ has a non-trivial $S^*$-symmetric infinitesimal motion  in $\H^{d}$.
\end{enumerate}
If $p$ is an $S$-regular point of $G_{H}$, then the following are equivalent:
\begin{enumerate}
\item[(i)]  $(G_{H},p)$ has an $S$-symmetric finite flex in $\E^{d}$;
\item[(ii)]  $(G_{H},q)$ has an $S^*$-symmetric finite flex in $\SS_{+}^{d}$;
\item[(iii)]  $(G_{H},r)$ has an $S^*$-symmetric finite flex in $\H^{d}$.
\end{enumerate}
\end{theorem}


\subsection{Inversion and the whole sphere}

When we consider the range of frameworks on the sphere, and of symmetries of frameworks on the sphere, it is not natural to restrict the joints to the upper half-sphere $\SS_{+}^{d}$.   We want to include joints which are on the lower half-sphere, and the symmetries which move points around on the entire sphere.

The processes of pulling and pushing orbits of joints on the cone include such points, provided we scale orbits of joints on the upper half-sphere by $\alpha=-1$.  An alternative interpretation  of this process is that we are applying the inversion $I(q) = - q$ to selected orbits.

It is clear from the analysis above that inverting an orbit of joints preserves the symmetries in the group $S^{*}$.  It has no impact on the rank of the corresponding cone orbit matrix, though the infinitesimal velocities will be multiplied by $-1$.  Applying this to a single orbit of vertices will also multiply the scalars of the self-stress for edges incident with the vertex representing the orbit in the cone orbit matrix by $-1$. We then have the following extension of Theorem~\ref{thm:hemisphere}:

\begin{cor}[Transfer between Euclidean and spherical spaces] \label{thm:wholesphere}
Given a configuration $q$ of $n$ points  in $\SS^{d}$ such that the projection $\pi(q)$ from the origin is equal to $p \in \E^{dn}$, we have that
\begin{enumerate}
\item[(i)]  the space of $S$-symmetric infinitesimal motions  of $(G,p)$ in $\E^{d}$ is isomorphic to the space of $S^{*}$-symmetric infinitesimal motions of $(G,q)$ in $\SS^{d}$;

\item[(ii)]  the space of $S$-symmetric trivial infinitesimal motions of $(G,p)$ in $\E^{d}$ is isomorphic to the space of $S^{*}$-symmetric trivial infinitesimal motions of $(G,q)$ in $\SS^{d}$;

\item[(iii)]   $(G,p)$ has a non-trivial $S$-symmetric  infinitesimal motion in $\E^{d}$  if and only if $(G,q)$  has a non-trivial $S^{*}$-symmetric  infinitesimal motion in $\SS^{d}$;

\item[(iv)]  the space of $S$-symmetric  self-stresses  of $(G,p)$ in $\E^{d}$ is isomorphic to the space of $S^{*}$-symmetric  self-stresses of $(G,q)$ in $\SS^{d}$.
\end{enumerate}
\end{cor}

While inverting entire orbits of joints cannot reduce the symmetry group of the framework, it can increase the symmetry group of the framework.  Consider the framework $(K_{4},q)$ realized as the vertices of a regular tetrahedron, with one face horizontal in the upper hemisphere.  This still projects to a framework $(K_{4},p)$ which has $C_{3}$ symmetry (there is a threefold rotational axis through the bottom joint and the midpoint of the horizontal face of $(K_{4},q)$).

The added symmetries of the tetrahedron do not survive the projection, as isometries of the projection.   However, they are still projective transformations within the projected space - and can play into the infinitesimal rigidity properties of the projected framework, since infinitesimal rigidity is a projective invariant \cite{CW, WWcones}.

When we consider finite motions of a coned framework, we recall the general property that pulling and pushing joints (or orbits of joints) along the cone rays leaves the finite motions unchanged.   This is true without reference to symmetry or any other analysis which originally predicted the finite motion of one of the cones (and therefore all of the cone-equivalent frameworks).

What is more surprising is that we could invert some, but not all joints of an orbit without any impact on the infinitesimal rigidity properties, or the finite flexibility.  Looked at in the other direction, when examining a framework on the sphere which may have no symmetries, we might still predict flexibility by noticing that the projection into the plane does have symmetry!   This says that key properties can be predicted by `seeing' the apparent symmetry of the structure when viewed from the center of the sphere!


\subsection{Other Cayley-Klein metrics}
\label{sec:projective}

The shared infinitesimal rigidity properties of  frameworks in $E^{d}$, $\SS^{d}$, and $\H^{d}$, as well as the connections to coning, all have a common geometric basis - a shared projective geometry.  For more than 150 years, there has been a recognition that infinitesimal rigidity is projectively invariant \cite{CW}.   In fact, the process of coning, applying an isometry, then reprojecting generates the projective transformations of a space (if repeated).

Given this shared projective geometry, it is natural to consider the other metrics which can be layered onto the projective geometry - the Cayley-Klein geometries \cite{CayleyKlein, richter}.  One approach to a number of these geometries is through quadratic forms.  In \S\ref{sec:Minkowski} we described the Minkowskian geometry as the metric with signature $(+,\ldots, +,-)$.  We then brought in the hyperbolic and de Sitter spaces as the geometry on `spheres' in this space.

There are analogous Minkowskian geometries $\M^{d+1}_{p,q}$ for other signatures $(p,q)=(+,\ldots, +, -,\ldots, -)$ (by convention  $p\geq q$). (Notice that we do not include cases, such as the Galilean geometry, with some $0$'s in the signature.)  These  geometries $\M^{d+1}_{p,q}$ have analogous rigidity matrices.  The hyperplane $x_{d+1}=1$ in $\M^{d+1}_{p,q}$ has the geometry of $\M^{d}_{p,q-1}$.   Each of these spaces has a distinct group of isometries - and therefore distinct point group symmetries.   However, the entire process used in \S\ref{sec:Minkowski} can be applied to show that the symmetric infinitesimal rigidity, and symmetry induced finite flexes at regular points, are equivalent in  the Cayley-Klein geometries $\M^{d}_{p,q-1}$ and in the corresponding Cayley-Klein geometries  on spheres in $\M^{d+1}_{p,q}$ \cite{CayleyKlein}.

\medskip

For simple infinitesimal flexibility, there is an additional transfer between the spaces  $\M^{d+1}_{p,q}$ and $\M^{d+1}_{p+k,q-k}$ for any $q-p \leq k \leq q$.  This is proven by  simply writing down the corresponding rigidity matrices and multiplying the corresponding $k$ columns by $-1$, which switches their signature by $-1$.   The rank of the matrices remains unchanged, and the kernels are isomorphic:  we simply multiply the corresponding entries of an infinitesimal motion by $-1$ to respect the change in signature.  This further emphasizes the underlying projective nature of the infinitesimal rigidity \cite{CW,SalWW}.   One overview of this connection involves the projective form of the rigidity matrix from which the transfer to these metrics is just a matter of taking the appropriate form of `perpendicular' for the metric in question, and writing a modified rigidity matrix of the same rank (see \cite{SalWW} for some background on this).

Transferring the symmetry results is more subtle - particularly the results for finite flexibility.   The first issue is whether the same symmetry group is even available in the corresponding spaces, as we move across signatures.  Two key examples of symmetries which always transfer are half-turn symmetry in $\M^{d+1}_{p,q}$ about an axis which is aligned with $d-1$ of the coordinate axes, or a reflection about a hyperplane with $x_{i}=0$. Notice that the transfer of half-turn symmetry is significant for predicting flexibility of structures such as one class of the Bricard octahedra in spaces of dimension $3$ \cite{alex,Stachelhyp,BSWWorbit}.   Similarly, mirror symmetry is key to predicting the flexibility of a second class of Bricard octahedra in these spaces of dimension $3$.    On the other hand, 3-fold rotation about the origin, for example, does not transfer from $\E^{2}$ to the Minkowskian plane $\M^{2}$.  There is a great deal to explore in these general transfers - !
 though the significance of rigidity in those spaces remains to be established.


\subsection{Tensegrity frameworks}

We first present a few basic definitions and some standard results for tensegrity frameworks in the symmetric setting.

A \emph{tensegrity graph} $\hat G$ has a partition of the edges of $G$ into three disjoint parts $E(G)=E_+(G)\cup E_-(G)\cup E_0(G)$.  $E_+(G)$ are the edges that are {\em cables}, $E_-(G)$ are the {\em struts} and $E_0(G)$ are the {\em bars}.  For a tensegrity framework $(\hat G, p)$, a {\em proper self-stress} is a self-stress on the underlying framework $(G,p)$ with the added condition that   $\omega_{ij}\geq 0, \{i,j\}\in E_+$, $\omega_{ij}\leq 0, \{i,j\}\in E_-$ \cite{RW1}.

Given a symmetric framework $(G,p)\in \mathscr{R}_{(G,S)}$, it is possible to use an $S$-symmetric self-stress on the bar-and-joint framework $(G,p)$ to investigate both the infinitesimal rigidity of $(G,p)$,  and  the infinitesimal rigidity of an associated \emph{$S$-symmetric} tensegrity framework $(\hat{G},p)$ (i.e., the edges of an edge orbit are either all cables, or all struts, or all bars),  with all members with $\omega_{ij}>0$ as cables and all members with $\omega_{ij}<0$ as struts \cite{BSWWorbit}.

The standard result for the infinitesimal rigidity of such tensegrity frameworks is:

\begin{theorem} [Roth, Whiteley \cite{RW1}]
\label{thm:RothWhiteley} A tensegrity framework $(\hat{G},p)$ is infinitesimally rigid if and only if the underlying bar framework $({G},p)$ is infinitesimally rigid as a bar-and-joint framework and $({G},p)$ has a self-stress which has  $\omega_{ij}>0$ on cables and  $\omega_{ij}<0$ on struts.
\end{theorem}

Translated in terms of the orbit matrix for a symmetric framework, this becomes \cite{BSWWorbit}:

\begin{cor}  \label{thm:SymmetricTensegrity} An  $S$-symmetric tensegrity framework $(\hat{G},p)$ is infinitesimally rigid if and only if the underlying bar framework $({G},p)\in \mathscr{R}_{(G,S)}$ is infinitesimally rigid as a bar-and-joint framework and the orbit matrix $\mathbf{O}({G},p,S)$ has an $S$-symmetric self-stress which has  $\omega_{ij}>0$ on cables and  $\omega_{ij}<0$ on struts.
\end{cor}

For this paper, what is of interest is the transfer through coning, with and without symmetry.  We have seen that the two conditions of Theorem~\ref{thm:RothWhiteley} and Corollary~\ref{thm:SymmetricTensegrity} are transferred by coning from $\E^{d}$ to $\E^{d+1}$, provided that we are pulling and pushing orbits of joints within the upper half-space:  (1) the infinitesimal rigidity of the underlying bar framework is preserved by the coning, and (2) the coefficients of the self-stress are multiplied by positive scalars so that the signs of the self-stress on orbits of original edges remain the same.
The signs of the self-stress (if any) on the edges to the cone joint are not as easily predicted, though they retain the symmetry of $S^{*}$.  In this simple coning of a tensegrity framework we will assume that the partition of the original edges into cables, struts, and bars are maintained, and that the edges to the cone joint become bars, creating $(\hat{G*o},q)$.

This  gives the following result:

\begin{theorem}[Tensegrity coning]  \label{thm:SymmetricTensegrityConing} Given an  $S$-symmetric tensegrity framework $(\hat{G},p)$, and a corresponding $S^{*}$-symmetric cone framework $(\hat{G*o},q)$  with $q$ in the upper half-space and $\pi(q)=p$, then  $(\hat{G},p)$ is infinitesimally rigid as a tensegrity framework in $\E^{d}$   if and only if
$(\hat{G*o},q)$ is  infinitesimally rigid as a tensegrity framework in $\E^{d+1}$.
\end{theorem}

Clearly, Theorem~\ref{thm:SymmetricTensegrityConing} has a corollary for the transfer to the upper hemisphere. We can also immediately extend Theorem~\ref{thm:SymmetricTensegrityConing} to cones into the upper half of Minkowskian space, giving a corollary in hyperbolic space.

\begin{cor}[Tensegrity transfer]  \label{thm:TensegrityTransfer} Given an  $S$-symmetric tensegrity framework $(\hat{G},p)$ in $\E^{d}$, and corresponding $S^{*}$-symmetric cone tensegrity frameworks $(\hat{G},q)$ and  $(\hat{G},r)$ in $\SS_{+}^{d}$ and $\H^{d}$, respectively,  with $q$ and $r$ in the upper half-space and $\pi(q)=p=\pi(r)$, then the following are equivalent:
\begin{enumerate}
\item[(i)] $(\hat{G},p)$ is infinitesimally rigid as a tensegrity framework in $\E^{d}$;

\item[(ii)] $(\hat{G},q)$ is infinitesimally rigid as a tensegrity framework in $\SS_{+}^{d}$;

\item[(iii)] $(\hat{G},r)$ is infinitesimally rigid as a tensegrity framework in $\H^{d}$.

\end{enumerate}
\end{cor}

If we want to work on  the whole sphere, using inversion of orbits as described above, we need to track the sign changes in the self-stress and make the corresponding changes in the partition of the edge orbits.  Specifically,
\begin{itemize}
\item if exactly one end-vertex of an edge is inverted, then there is a change in the sign of the stress and  we switch the assignment of the edge orbit in the partition of cables or struts;
\item if both end-vertices of an edge are inverted, then there is no net change in the sign of the stress and we make no change in the assignment of the edge orbit in the partition;
\item if neither end-vertex of an edge is inverted, then there is no change in the sign of the stress and we make no change in the assignment of the edge orbit in the partition.
\end{itemize}

While there is more that can be said, we hope this brief summary provides the flavor of how coning impacts the rigidity properties of tensegrity frameworks.


\providecommand{\bysame}{\leavevmode\hbox to3em{\hrulefill}\thinspace}
\providecommand{\MR}{\relax\ifhmode\unskip\space\fi MR }
\providecommand{\MRhref}[2]{%
  \href{http://www.ams.org/mathscinet-getitem?mr=#1}{#2}
}
\providecommand{\href}[2]{#2}

\end{document}